\documentclass[letterpaper,12pt]{article}
\usepackage{amsfonts,latexsym,amssymb,amsmath,algorithmic,epsfig}
\newlength{\jmr}
\newlength{\barvinok}
\newtheorem{thm}{Theorem}
\newtheorem{lemma}{Lemma}
\newtheorem{cor}{Corollary}
\newtheorem{conj}{Univariate Threshold Conjecture.}

\newtheorem{ks}{Karpinski-Shparlinski Theorem}

\renewcommand{\mod}{\mathrm{mod}}
\newtheorem{baker}{Baker's Theorem}

\newtheorem{algsimp}{Algorithm {\tt DimPlus1Case}}

\newtheorem{algcirc}{Algorithm {\tt DimPlus2Case}}

\newtheorem{algbinv}{Algorithm {\tt BinomialVanish}}

\newtheorem{algbins}{Algorithm {\tt BinomialSign}}

\newtheorem{rems}{Remark} 
\newtheorem{dfn}{Definition}

\newtheorem{ex}{Example}
\newtheorem{prop}{Proposition}

\newlength{\smale}
\newcommand{\thth}{^{\text{\underline{th}}}}

\newcommand{\np}{{\mathbf{NP}}}

\newcommand{\rp}{{\mathbf{RP}}}
\newcommand{\feas}{{\mathbf{FEAS}}}
\newcommand{\disc}{{\mathbf{ADISCVAN}}}
\newcommand{\sdisc}{{\mathbf{ADISCSIGN}}}

\newcommand{\gln}{\mathbb{G}\mathbb{L}_n}
\newcommand{\glm}{\mathbb{G}\mathbb{L}_m}

\newcommand{\bpp}{{\mathbf{BPP}}}

\newcommand{\vol}{{\mathrm{Vol}}}

\newcommand{\sat}{\mathbf{3CNFSAT}} 
\newcommand{\pp}{\mathbf{P}}

\newcommand{\nc}{\mathbf{NC}}
\newcommand{\pspa}{\mathbf{PSPACE}}

\newcommand{\eps}{\varepsilon}
\newcommand{\Pro}{{\mathbb{P}}}

\newcommand{\barf}{{\bar{f}}}

\newcommand{\Q}{\mathbb{Q}}
\newcommand{\R}{\mathbb{R}}
\newcommand{\C}{\mathbb{C}}
\newcommand{\N}{\mathbb{N}}
\newcommand{\Z}{\mathbb{Z}}
\newcommand{\bO}{\mathbf{O}}

\newcommand{\Zn}{\Z^n}

\newcommand{\Rn}{\R^n}

\newcommand{\Cs}{\C^*}
\newcommand{\Rs}{\R^*}
\newcommand{\cA}{{\mathcal{A}}}
\newcommand{\cB}{{\mathcal{B}}}
\newcommand{\cC}{{\mathcal{C}}}

\newcommand{\cE}{{\mathcal{E}}}
\newcommand{\cF}{{\mathcal{F}}}

\newcommand{\cM}{{\mathcal{M}}}
\newcommand{\cN}{{\mathcal{N}}}

\newcommand{\cNc}{\cN_{\mathrm{comp}}}
\newcommand{\cNn}{\cN_{\mathrm{non}}}
\newcommand{\cQ}{{\mathcal{Q}}}
\newcommand{\cQS}{{\cQ\cS}}
\newcommand{\cS}{{\mathcal{S}}}

\newcommand{\cU}{{\mathcal{U}}}
\newcommand{\cV}{{\mathcal{V}}}

\newcommand{\Rsn}{{(\R^*)}^n}

\newcommand{\Csn}{{(\C^*)}^n}
\newcommand{\qed}{$\blacksquare$}
\newcommand{\dia}{$\diamond$}

\newcommand{\newt}{\mathrm{Newt}}
\newcommand{\codim}{\mathrm{codim}}
\newcommand{\size}{\mathrm{size}}
\newcommand{\conv}{\mathrm{Conv}}

\newcommand{\supp}{\mathrm{Supp}}
\newcommand{\init}{\mathrm{In}}
\newcommand{\sign}{\mathrm{sign}}

\begin{document}
\title{\mbox{}\\
\vspace{-1in}First Steps in Algorithmic Real Fewnomial Theory } 
\author{
Frederic Bihan\thanks{UFR SFA, Campus Scientifique, 73376 Le Bourget-du-Lac 
Cedex, France. {\tt Frederic.Bihan@univ-savoie.fr}  , 
{\tt www.lama.univ-savoie.fr/\~{}bihan} \ .  } 
\and 
J. Maurice Rojas\thanks{Department of Mathematics,
Texas A\&M University
TAMU 3368,
College Station, Texas 77843-3368,
USA. {\tt rojas@math.tamu.edu } , {\tt www.math.tamu.edu/\~{}rojas} \ .
Partially supported by 
NSF individual grant DMS-0211458, NSF CAREER grant 
DMS-0349309, Sandia National Laboratories, and the American 
Institute of Mathematics. } \and Casey E.\ 
Stella\thanks{Partially supported by NSF grant DMS-0211458.}} 

\date{\today} 

\maketitle

\vspace{-.5cm}
\noindent
\mbox{}\hfill
{\scriptsize J.\ Maurice Rojas dedicates this paper to the memory of his dear 
friend, Richard Adolph Snavely, 1955--2005.}\hfill
\mbox{}\\

\vspace{-.5cm}
\begin{abstract} 
Fewnomial theory began with explicit bounds --- solely in terms of 
the number of variables and monomial terms --- on the number of 
real roots of systems of polynomial equations.  
Here we take the next logical step of investigating the corresponding
existence problem: Let $\feas_\R$ denote the problem of deciding whether a
given system of multivariate polynomial equations with integer coefficients 
has a real root or not. We describe a phase-transition for when $m$ is 
large enough to make $\feas_\R$ be $\np$-hard, when restricted to 
inputs consisting of a single $n$-variate polynomial with exactly $m$ monomial 
terms: polynomial-time for $m\!\leq\!n+2$ (for any fixed $n$) and 
$\np$-hardness for $m\!\geq\!n+n^\eps$ (for $n$ varying and any 
fixed $\eps\!>\!0$). Because of important connections 
between $\feas_\R$ and $A$-discriminants, we then study some new families of 
$A$-discriminants whose signs 
can be decided within polynomial-time. ($A$-discriminants contain all known 
resultants as special cases, and the latter objects are central in algorithmic 
algebraic geometry.) Baker's Theorem from diophantine approximation 
arises as a key tool. Along the way, we also derive new quantitative
bounds on the real zero sets of $n$-variate $(n+2)$-nomials.
\end{abstract} 

\section{Introduction and Main Results} 
Let $\feas_\R$ --- a.k.a.\ the {\bf real feasibility problem} --- 
denote the problem of deciding whether a given system of 
polynomial equations with integer coefficients has a real root or not.
While $\feas_\R$ is arguably the most fundamental problem of real algebraic
geometry, our current knowledge of its computational complexity is surprisingly
coarse, especially for sparse polynomials. This is a pity, for in addition to 
numerous practical applications (see, e.g., \cite{dimacs}), $\feas_\R$ 
is also an important motivation behind effectivity estimates for the Real 
Nullstellensatz (e.g., \cite{stengle,schmid}), the quantitative study of 
sums of squares \cite{blekherman}, and their connection to semi-definite 
programming \cite{parrilo}. Furthermore, efficient algorithms for $\feas_\R$ 
are crucial for the tractability of harder problems such as 
quantifier elimination, and computing the closure and frontier of 
more general types of varieties such as sub-Pfaffian sets \cite{gv}.

Before stating our main results, we will need to clarify some 
geometric notions concerning sparse polynomials. 
\begin{dfn}
Let $f(x)\!:=\!\sum^m_{j=1} c_jx^{a_j}\!\in\!\Z[x^{\pm 1}_1,\ldots,
x^{\pm 1}_n]$, where $x^{a_j}\!:=\!x^{a_{1j}}_1\cdots x^{a_{nj}}_n$, 
$c_j\!\neq\!0$ for all $j$, and the $a_j$ are distinct. We call such an $f$ 
an {\bf $\pmb{n}$-variate $\pmb{m}$-nomial} and
we call $\{a_1,\ldots,a_m\}$ the {\bf support} of $f$.
Also, for any collection $\cF$ of polynomial systems with\linebreak
\scalebox{.96}[1]{integer coefficients, let $\feas_\R(\cF)$ denote the natural
restriction of $\feas_\R$ to inputs in $\cF$. \dia} 
\end{dfn}
Note that $n$-variate quadratic polynomials are a special case of 
$n$-variate $m$-nomials with $m\!=\!O(n^2)$. 
\begin{dfn} 
\scalebox{.87}[1]{Let us use $\#S$ for the cardinality of a set $S$ and say 
that a subset $A\!=\!\{a_1,\ldots,a_m\}\!\subset\!\Zn$}\linebreak  
\scalebox{.98}[1]{with $\#A\!=\!m$ is {\bf affinely independent} iff  
the implication ``$[\sum^m_{j=1} \gamma_j a_j=\bO$ and $\sum^m_{j=1} 
\gamma_j\!=\!0]$}\linebreak $\Longrightarrow$ 
$\gamma_1\!=\cdots=\!\gamma_m\!=\!0$''
holds for all $(\gamma_1,\ldots,\gamma_m)\!\in\!\R^m$. 
Also let $\dim A$ denote   
the dimension of the subspace of $\Rn$ 
generated by the set of all {\bf differences} of vectors in $A$. \dia 
\end{dfn}

\noindent
Note in particular that $\#A\!\geq\!1+\dim A$, with  
equality iff $A$ is affinely independent. It is also easily checked that a 
random set of $m$ points in $\Rn$ (chosen, say, independently from any 
continuous probability distribution) will have dimension $\min\{n,m-1\}$ 
with probability $1$. 

Clearly, any $n$-variate Laurent polynomial $f$ always 
satisfies $\dim \supp(f)\!\leq\!n$.  Let $\Rs\!:=\!\R\setminus\{0\}$. 
It is not much harder to see that one can always find (even in an 
algorithmically efficient sense) a $d$-variate Laurent polynomial 
$g$, with the same number of terms as $f$ and $d\!=\!\dim\supp(f)$, such that 
$f$ vanishes in $\Rsn$ iff $g$ vanishes in  $(\Rs)^d$ (see Corollary 
\ref{cor:hyp} of Section \ref{sec:back} below). In this sense, 
``almost all'' $n$-variate $m$-nomials satisfy 
$\dim \supp(f)\!=\!\min\{n,m-1\}$, and those that don't are essentially 
just $d$-variate $m$-nomials (with $d\!<\!n$) in disguise.  

Recall the containments of complexity classes
$\nc_1\!\subseteq\!\nc\!\subseteq\!\pp\!\subseteq\!\rp\!\subseteq\!
\bpp\cup\np\!\subseteq\!
\pspa$ (these complexity classes are reviewed briefly in Section \ref{sec:back} 
below). 
Roughly speaking, our first main result  says that real feasibility for 
``honest'' $n$-variate $(n+k)$-nomials is easy for $k\!\leq\!2$, but 
$\np$-hardness kicks in quickly already for $k$ a very 
slowly growing function of $n$.
\begin{thm}
\label{thm:me} 
Let 

\vspace{-.3cm}
\noindent
\begin{eqnarray*}
\cA & := & \left.  \{f\!\in\!\Z[x_1,\ldots,x_n] \; \right| \; 
\#\supp(f)\!=\!1+\dim\supp(f) \text{ and } n\!\in\!\N\},\\
\cB_n & := & \left. \{f\!\in\!\Z[x_1,\ldots,x_n] \; \right| \; 
\#\supp(f)\!=\!2+\dim\supp(f)\}, \\
\cC_\eps & := & \scalebox{1}[1]{\text{$\left\{f\!\in\!\Z[x_1,\ldots,x_n] 
\; \left| \; \begin{matrix}\#\supp(f)\!\leq\!n+n^\eps,  \ \dim\supp(f)\!=\!n,  
\ n\!\in\!\N, \\
\text{and } f \text{ is a sum of squares of polynomials.} \end{matrix}
\right\}\right., \text{ and}$}} \\
\cS & := & \text{\scalebox{1}[1]{$\left\{(f_1,\ldots,f_k) \; \left| \; 
\begin{matrix} f_i\!\in\!\Z[x_1,\ldots,x_n] 
\text{ and } f_i \text{ is a linear trinomial or a}\\
\text{binomial of degree } \leq\!2 \text{ for all } i, 
\text{ and } k\!\geq\!n\!\geq\!1.\end{matrix}\right. \right\}$}}  
\end{eqnarray*}

\vspace{-.3cm}
\noindent
Then, measuring the {\bf size} of any polynomial 
$f(x)\!=\!\sum^m_{j=1}c_jx^{a_j}\!\in\!\Z[x_1,\ldots,x_n]$ --- 
denoted $\pmb{\size(f)}$ --- as 
the total number of binary digits in the $c_i$ and $a_{i,j}$, we have:\\
\vspace{-.5cm} 
\begin{enumerate}
\item{ $\feas_\R(\cA)\!\in\!\nc_1$. }
\vspace{-.3cm} 
\item{For any {\bf fixed} $n\!\in\!\N$, $\feas_\R(\cB_n)\!\in\!\pp$. } 
\vspace{-.3cm} 
\item{$\feas_\R(\cS)$ is $\np$-hard.} 
\vspace{-.3cm} 
\item{ For any {\bf fixed} $\eps\!>\!0$, $\feas_\R(\cC_\eps)$ is $\np$-hard. } 
\end{enumerate}
\end{thm}
\begin{ex}
A very special case of Assertion (2) of Theorem \ref{thm:me} implies that 
one can decide --- for any nonzero $c_1,\ldots,c_5\!\in\!\Z$ and 
$D\!\in\!\N$ --- whether \\
\mbox{}\hfill$c_1+c_2x^{999}+c_3x^{73}z+c_4y^D+c_5x^Dy^{3D}z^{9D}$
\hfill\mbox{}\\
\scalebox{.91}[1]{has a root in $\R^3$, using a number of bit operations 
polynomial in $\log(D)+\log\left[(|c_1|+1)\cdots (|c_5|+1)\right]$.}\linebreak 
The best previous results (e.g., via the critical points method, 
infinitesimals, and rational univariate reduction, as detailed in 
\cite{bpr}) would yield a bound polynomial in\linebreak  
$D+\log\left[(|c_1|+1)\cdots (|c_5|+1)\right]$ instead. Assertion (2) 
also vastly generalizes an earlier analogous result for univariate trinomials 
\cite{rojasye}. \dia
\end{ex}

\noindent 
The algorithm underlying Assertion (2) turns out to depend critically on the 
combinatorics of $\supp(f)$, particularly its triangulations. Furthermore, 
extending the polynomiality of 
$\feas_\R(\cF)$ from $(n+1)$-nomials to $(n+2)$-nomials turns out to be 
surprisingly intricate, involving $A$-discriminants (cf.\ Section 
\ref{sub:disc}), Baker's Theorem on Linear Forms in Logarithms (cf.\ 
Section \ref{sec:proofs}), and Viro's Theorem from toric geometry 
(see, e.g., \cite[Thm.\ 5.6]{gkz94}). Theorem 1, along with some more 
technical strengthenings, is proved in Section \ref{sub:proof1} below. 
\begin{rems} 
There appears to have been no 
earlier explicit statement that $\feas_\R(\cF)\!\in\!\pp$ (or 
even $\feas_\R(\cF)\!\in\!\np$) for $\cF$ some 
non-trivial family of $n$-variate $m$-nomials with $m\!=\!n+O(1)$. 
As for lower bounds, the best previous result for sparse polynomials 
appears to have been $\np$-hardness of $\feas_\R(\cF)$ 
when $\cF$ is the family of those $n$-variate $m$-nomials with 
$m\!=\!\Omega(n^3)$ (see, e.g., \cite[discussion preceding 
Thm.\ 2]{rojasye}). Assertion (4) is therefore a considerable 
sharpening.\footnote{An earlier version of this paper 
proved $\np$-hardness of $\feas_\R(\cF)$ when $\cF$ is the family of 
$n$-variate $(6n+6)$-nomials, but an anonymous referee suggested an even 
easier proof for the stronger version of Assertion (4) we are now stating.} 
\dia 
\end{rems}  
\begin{rems} 
Let $\cU_m\!:=\!\{f\!\in\Z[x_1] \; | \; f \text{ has exactly } m  
\text{ monomial terms}\}$. 
While it has been known since the late 1980's that $\feas_\R\!\in\!\pspa$
\cite{pspace}, it is already unknown whether
$\feas_\R(\cU_4)\!\in\!\bpp\cup\np$, or even whether 
$\feas_\R(\Z[x_1,\ldots,x_n])$ is $\np$-hard for some particular 
value of $n$ \cite{lickroy,rojasye}. (The latter reference nevertheless 
states certain analytic hypotheses under which it would follow that 
$\feas_\R(\cU_m)\!\in\!\pp$ for fixed $m$.) 
The role of sparsity in complexity bounds for univariate real feasibility  
is thus already far from trivial.  \dia 
\end{rems}  

Quantitative results over $\Rn$ of course shape the kind 
of algorithms we can find over $\Rn$. In particular,  
Khovanski's famous {\bf Theorem on Real Fewnomials} implies an upper bound
depending only on $m$ and $n$ --- {\bf independent} of the degree ---  
for the number of connected components of the real zero set of
any $n$-variate $m$-nomial \cite{few}.
More recently, his bound has been improved from $2^{O(m^2)}n^{O(n)}n^{O(m)}$ 
in the smooth case \cite[Sec.\ 3.14, Cor.\ 5]{few} 
to $2^{O(m^2)}2^{O(n)}n^{O(m)}$ in complete generality
\cite[Cor.\ 2]{tri} (see also \cite{perrucci} for further improvements). 

For $n$-variate $(n+2)$-nomials we can now make a dramatic improvement. 
Recall that a set $S\!\subseteq\!\Rn$ is {\bf convex} iff
for any $x,y\!\in\!S$, the line segment connecting $x$ and $y$
is also contained in $S$. Recall also that for any $A\!\subset\!\Rn$, the 
{\bf convex hull of $A$} --- denoted $\conv A$ --- is the smallest convex set
containing $A$.
\begin{thm}
\label{thm:ckt} 
Let $f$ be any $n$-variate $m$-nomial with $m\!\leq\!2+\dim\supp(f)$, 
$Z_+(f)$ its zero set in $\Rn_+$, and define $\cNc(f)$ (resp.\ 
$\cNn(f)$) to be the number of compact (resp.\ non-compact) connected 
components of $Z_+(f)$. Then $\cNc(f)\!\leq\!1$ (with  
examples attaining equality for each $n\!\in\!\N$), and  
$\cNc(f)\!=\!1 \Longrightarrow Z_+(f)$ is either 
a point, or isotopic to an $(n-1)$-sphere. 
Also, $\cNn(f)$ is no more than $0$, $2$, $6$, $9$, or 
$2n+2$, according as $n$ is $1$, $2$, $3$, $4$, or $\geq\!5$ (with 
examples attaining equality for each $n\!\leq\!2$). 
Finally, if $\supp(f)$ is disjoint from the interior of its convex 
hull, then $\cNc(f)\!=\!0$.  
\end{thm} 
While the bound $\cNn(f)\!\leq\!2$ (for $n\!=\!2$) was found earlier by 
Daniel Perrucci, all the other bounds of Theorem \ref{thm:ckt} are new, and 
the special case $n\!=\!3$ improves Perrucci's earlier bound for 
$3$-variate $5$-nomials by a factor of at least $2$ 
(see \cite[Thms.\ 4 \& 5]{perrucci}). Except for an upper 
bound of $2^{O(n^2)}$ for the 
smooth case \cite[Sec.\ 3.14, Cor.\ 5]{few}, there appear to be no other 
earlier explicit bounds in the spirit of Theorem \ref{thm:ckt}.
Theorem \ref{thm:ckt} has recently been generalized to systems of $k$ 
polynomials in $\cB_n$, with identical supports 
and $k\!\in\!\{1,\ldots,n\}$ \cite{brs}. (See also \cite{bbs,bihan} 
for the opposite extreme to Theorem \ref{thm:ckt}: bounding the number of 
isolated real roots of $n$ polynomials in $\cB_n$ with the same support.) 
We prove Theorem \ref{thm:ckt} in Section \ref{sub:proof2}. 

It thus appears that, unlike algebraic geometry over $\C$, 
large degree is potentially less of a complexity bottleneck over $\R$. 
Considering the ubiquity of sparse real polynomial systems in engineering, 
algorithmic speed-ups in broader generality via sparsity are thus of  
the utmost interest. Furthermore, in view of the complexity 
threshold of Theorem \ref{thm:me}, randomization, approximation, 
and/or average-case speed-ups   
appear to be the next key steps if we are to have a sufficiently general 
and useful algorithmic fewnomial theory over $\R$. A promising 
step in this direction can be found in work of Barvinok \cite{approx}, 
but more work still needs to be done before we can assert significant 
new randomized algorithms --- even for $\feas_\R(\cU_4)$. 

Since the largest $m$ for which $\feas_\R$ is doable in polynomial-time  
--- for input a single $n$-variate $m$-nomial --- appears to be 
$m\!=\!n+2$ (as of early 2006), we propose the following conjecture to address 
the cases $m\!\geq\!n+3$.
\begin{conj}
For any $m\!\in\!\N$, let\\
\mbox{}\hfill  $\cU_m\!:=\!\{f\!\in\!\Z[x_1] \; | \; f \text{ has exactly } m 
\text{ monomial terms}\}$.\hfill\mbox{}\\ 
Then $\feas_\R(\Z[x_1])$ is $\np$-hard but, for any {\bf fixed} $m$, 
there is a natural probability measure on $\cU_m$ so that 
$\feas_\R(\cU_m)$ has polynomial-time complexity on average. 
\end{conj} 

\noindent 
Note that $\feas_\R(\cU_3)\!\in\!\pp$, thanks to Assertion (2) of  
Theorem \ref{thm:me}, since $U_3\!=\!\cB_1$. The latter 
``positive'' part of the conjecture is meant to be reminiscent of 
Smale's 17$\thth$ Problem, which concerns the complexity of approximating 
complex roots of polynomial systems \cite{21,rojasye}. 
 
That feasibility over $\R$ may be $\np$-hard already for {\bf uni}variate 
sparse polynomials is suggested by a recent parallel over a different 
complete field: $\Q_p$. In particular, (a) there is 
now a {\bf Theorem on $p$-adic Fewnomials} (due to the middle author 
\cite{amd}), with significantly sharper bounds than Khovanski's Theorem on 
Real Fewnomials, and (b) it is now known that 
$\feas_{\Q_p}(\Z[x_1])$ --- the natural $p$-adic analogue of 
$\feas_{\R}(\Z[x_1])$ --- is doable in randomized polynomial-time only if 
an unlikely containment of complexity classes occurs: $\np\!\subseteq\!\bpp$ 
\cite{pr}. 

Possible alternative evidence for the Univariate Threshold 
Conjecture can be given via the {\bf $A$-discriminant}, which 
is defined in Section \ref{sub:disc} below. To set the stage, 
first recall the following elementary example.
\begin{ex} 
Note that $f(x_1)\!:=\!a+bx_1+cx^2_1$ has either $0$ or $2$ real roots 
according as the discriminant $\Delta\!=\!b^2-4ac$ 
is negative or positive. Observe then that the real zero set, 
$\widetilde{W}\!\subset\!\R^3$, of $\Delta$ can be identified with the 
collection of all quadratic polynomials possessing a degenerate root, and 
that $\Delta$ also defines a curve $W$ in the real projective plane 
$\Pro^2_\R$. Furthermore, $W$ is 
equivalent (under a linear change of variables over $\Q$) to a circle, 
and thus there are exactly $2$ discriminant chambers. These $2$ chambers 
correspond exactly to those quadratic polynomials possessing either $0$ or 
$2$ real roots. Finally, note that the support of $f$ is $\{0,1,2\}$ and this 
set admits exactly $2$ triangulations with vertices in $\{0,1,2\}$: They are 
\raisebox{-.05cm}{\epsfig{file=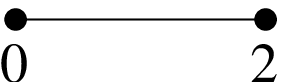,width=.5in}} and 
\raisebox{-.05cm}{\epsfig{file=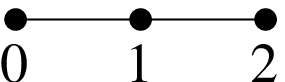,width=.5in}} . \dia  
\end{ex} 

\subsection{Complexity and Topology of Certain $A$-Discriminants} 
\label{sub:disc} 
The connection between computational complexity (e.g., of 
deciding membership in semi-algebraic sets \cite{dl79} 
or approximating the roots of univariate polynomials \cite{smaletopo}) 
and the topology of discriminant complements dates back to the late 1970's.  
Here, we point to the possibility of a more refined connection between 
$\feas_\R$ and discriminant complements. (See also \cite{drrs} for 
further results in this direction.) In particular, our last example was 
a special case of a much more general invariant attached to spaces of sparse 
multivariate polynomials. 
\begin{dfn} 
\label{dfn:adisc} 
Given any $A\!=\!\{a_1,\ldots,a_m\}\!\subset\!\Zn$ of 
cardinality $m$, define the set\\ 
\mbox{}\hfill
\scalebox{1}[1]{$\nabla^0_A\!:=\!\left\{(c_1,\ldots,c_m)\!\in\!\C^m\; 
| \; f(x)\!:=\!\sum^m_{j=1} c_j x^{a_j} \text{ has a  
degenerate root$^3$ in } \Csn 
\right\}$.}\setcounter{footnote}{3}\footnotetext{\scalebox{.9}[1]{That is, a 
root $\zeta$ of $f$ with $\frac{\partial f}
{\partial x_1}|_{x=\zeta}\!=\cdots=\!
\frac{\partial f}{\partial x_n}|_{x=\zeta}\!=\!0$.}}\hfill\mbox{}\\
\scalebox{.92}[1]{The {\bf $\pmb{A}$-discriminant} is then the unique 
(up to sign) irreducible polynomial 
$\Delta_A\!\in\!\Z[c_1,\ldots,c_m]\!\setminus\!\{0\}$}\linebreak whose 
complex zero set contains $\nabla^0_A$. (If $\codim \nabla_A\!>\!1$ then 
we set $\Delta_A\!:=\!1$.) For \linebreak convenience, we will usually write 
$\Delta_A(f)$ in place of $\Delta_A(c_1,\ldots,c_m)$. Finally, 
we let $\nabla_A$ --- the {\bf $\pmb{A}$-discriminant variety} --- 
denote the zero set of $\Delta_A$ in $\C^m$. \dia 
\end{dfn} 
\begin{ex} 
\scalebox{1}[1]{If we take $A\!:=\!\left\{
\text{\scalebox{1}[.8]{$\begin{bmatrix}0\\0\end{bmatrix},\ldots,\begin{bmatrix}d\\0\end{bmatrix},
\begin{bmatrix}0\\1\end{bmatrix},\ldots,\begin{bmatrix}e\\1\end{bmatrix}
$}}\right\}$}  then $\Delta_A(a_0,\ldots,a_d,b_0,\ldots,b_e)$ 
is exactly the classical 
{\bf Sylvester resultant} of the univariate polynomials $a_0+
\cdots+a_dx^d$ and $b_0+\cdots+b_ex^e$. This is a special 
case of a more general construction which shows that any multivariate 
{\bf toric resultant} can be obtained as a suitable 
$A$-discriminant \cite[The Cayley Trick, Prop.\ 1.7, pp.\ 274]{gkz94}. \dia 
\end{ex} 

\begin{dfn}
Let $\cF_A\!:=\!\{f\!\in\!\Z[x_1,\ldots,x_n]\; | \; \supp(f)\!\subseteq\!A\}$, 
and let $\disc$ (resp.\ $\sdisc$) denote the problem of deciding whether
$\Delta_A(f)$ vanishes (resp.\ determining\linebreak 
\scalebox{.94}[1]{the sign of $\Delta_A(f)$) for an 
input $f\!\in\!\cF_A$. Finally, let $\disc(\cF)$ (resp.\ $\sdisc(\cF)$)} 
\linebreak 
be the natural restriction of $\disc$ (resp.\ $\sdisc$) to input polynomials 
in some family $\cF$. \dia 
\end{dfn} 
\scalebox{.9}[1]{An intriguing link between $\feas_\R$ and 
discriminants is the fact that those $A$ with 
$\feas_\R(\cF_A)\!\in\!\pp$}\linebreak currently appear to coincide with 
those $A$ with $\disc(\cF_A)\!\in\!\pp$. 
This should not be too surprising in view of the following fact: 
If $f\!\in\!\cF_A$ has smooth complex zero set (and a similarly mild 
condition holds for its zero set at infinity), then $f$ lies in 
some connected component $C$ of $(\Rs)^{\#A}\setminus\nabla_A$ (under a 
natural identification of coefficients of $f$ and coordinates of 
$(\Rs)^{\#A}$), and any other $g\!\in\!C$ has real zero set 
isotopic to that of $f$ (see, e.g., \cite[Ch.\ 11, Sec.\ 5A, Prop.\ 5.2, 
pg.\ 382]{gkz94}). Let us call any such $C$ a {\bf discriminant chamber 
(for $A$)}. 

So deciding $\feas_\R(\cF_A)$ 
for a given $f$ is nearly the same as deciding whether $f$ lies 
in a particular union of discriminant chambers, and 
it is thus natural to suspect that the following 
three situations may be equivalent in some rigourous and useful sense: 
(a) $(\Rs)^{\#A}\setminus\nabla_A$ has ``few'' connected components, 
(b) $\feas_\R(\cF_A)$ is ``easy'', (c) $\disc(\cF_A)$ is ``easy''. 

While connections between (a) and (c) are known (see, e.g., 
\cite[Ch.\ 16]{bcss}), the best current theorems appear to be 
too weak to yield any complexity lower bounds of use for the 
Univariate Threshold Conjecture. As for stronger connections 
between (b) and (c), concrete examples arise, for instance, from the 
following two facts: (1) $\feas_\R(\cQ)\!\in\!\pp$, where 
$\cQ\!:=\!\{f\!\in\!\Z[x_1,\ldots,x_n] \; | \; f \text{ is homogeneous and 
quadratic}\}$ \cite{barvinok}\footnote{Barvinok actually proved 
the stronger fact that $\feas_\R(\cQS_k)\!\in\!\pp$ for any 
fixed $k$, where $\cQ\cS_k$ is the family of polynomial systems of the 
form $(f_1,\ldots,f_k)$ with $f_i\!\in\!\cQ$ for all $i$.} 
and (2) for $A$ the support of a quadratic polynomial $f$, $\Delta_A(f)$ 
is computable in polynomial-time.  (The latter fact follows easily from an 
exercise in Cramer's Rule for linear equations and the Newton identities  
\cite[Ch.\ 15, Pgs.\ 292--296]{bcss}.)   

\scalebox{.93}[1]{A new connection we can assert between the ``easiness'' of 
$\feas_\R(\cF_A)$ and $\disc(\cF_A)$}\linebreak is the case of $n$-variate 
$(n+2)$-nomials, thanks to Assertion (2) of Theorem \ref{thm:me} above, and 
the first part of our final main result below.   
\begin{thm} 
\label{thm:disc} 
Following the notation of Theorem \ref{thm:me} and our last 
two definitions: 
\begin{enumerate}
\item{$\disc\left(\bigcup\limits^\infty_{n=1}\cB_n\right)\!\in\!\pp$. } 
\item{For any {\bf fixed} $n$, $\sdisc(\cB_n)\!\in\!\pp$. }
\end{enumerate} 
\end{thm} 

\noindent 
Note that Theorem \ref{thm:disc} improves considerably on 
what can be done through quantifier elimination (e.g., \cite{pspace,bpr}), 
because for fixed $n$ these older methods already have complexity exponential 
in our notion of input size.  
Theorem \ref{thm:disc} --- proved in Section \ref{sec:disc} --- turns out to 
be a central tool in the algorithms behind the complexity upper bounds of 
Theorem \ref{thm:me}, and is the main reason that diophantine approximation 
enters our scenery. 

In light of the connections between the 
complexity of $\feas_\R(\cF_A)$ and $\disc(\cF_A)$, we 
conclude our introduction with another possible piece of evidence in favor 
of the ``negative'' portion of the Univariate Threshold Conjecture: 
\begin{ks}
\cite{ks} 
$\disc(\Z[x_1])$ is computationally hard in the following sense: If 
$\disc(\Z[x_1])\!\subseteq\!\cC$ for some 
complexity class $\cC$, then $\np\!\subseteq\!\cC\cup\rp$.\footnote{The 
paper \cite{ks} actually asserts the stronger fact that 
$\disc(\Z[x_1])$ is $\np$-hard, but without a proof. One of the authors of 
\cite{ks} (Igor Shparlinski) has confirmed this oversight, along with the 
fact that it was also observed independently by Erich Kaltofen 
\cite{shpar}.} In particular, $\disc(\Z[x_1])\!\subseteq\!\bpp 
\Longrightarrow \np\!\subseteq\!\bpp$. \qed 
\end{ks}

\noindent 
\scalebox{.86}[1]{The containment $\np\!\subseteq\!\bpp$ is widely 
disbelieved, so it would appear possible that $\disc(\Z[x_1])$}\linebreak 
is {\bf not} doable in randomized polynomial-time. 

Our main results are proved in Section \ref{sec:proofs}, after
the development of some necessary theory in Section \ref{sec:back} below. 
A useful elementary result on the real zero sets of $n$-variate 
$(n+1)$-nomials is then proved in Section \ref{sec:loose}. 
 
\section{Background and Ancillary Results} 
\label{sec:back}
Let us first informally review some well-known complexity classes (see, e.g., 
\cite{papa} for a complete and rigourous description). 
\begin{rems} 
Throughout this paper, our algorithms will always have a notion 
of input size that is clear from the context, and our underlying 
computational model will always be the classical Turing model \cite{papa}. 
Thus, appellations such as ``polynomial-time'' are to be understood 
as ``having bit-complexity polynomial in the underlying input size'', 
and the underlying polynomial and/or $O$-constants depend only 
on the algorithm, not on the specific\linebreak \scalebox{.97}[1]{instance 
being solved. The same of course applies to ``linear-time'', 
``exponential-time'', etc. \dia}  
\end{rems} 
\begin{itemize}
\item[$\nc_1$]{ The family of decision problems which can be done within 
time $O(\log \text{{\tt InputSize}})$, using a number of 
processors linear in the input size.  }
\item[$\nc$]{ The family of decision problems which can be done within time
poly-logarithmic in the input size, using a number of processors 
polynomial in the input size. }
\item[$\pp$]{ The family of decision problems which can be done within 
polynomial-time. }
\item[$\rp$]{ The family of decision problems admitting 
randomized polynomial-time algorithms 
for which a {\tt ``Yes''} answer is always correct but a {\tt ``No''} answer
is wrong with probability $\frac{1}{2}$. }
\item[$\bpp$]{ The family of decision problems admitting randomized 
polynomial-time algorithms that terminate with an answer that is correct with
probability at least\footnote{It is easily shown that we can
replace $\frac{2}{3}$ by any constant strictly greater than $\frac{1}{2}$
and still obtain the same family of problems \cite{papa}.} $\frac{2}{3}$.}
\item[$\np$]{ The family of decision problems where a {\tt ``Yes''} answer can
be {\bf certified} within polynomial-time.}
\item[$\pspa$]{ The family of decision problems solvable within 
polynomial-time, provided a number of processors exponential in
the input size is allowed. }
\end{itemize} 

\noindent
Recall also that even the containment $\pp\!\subseteq\!\pspa$ is still an open 
problem (as of early 2006).  

A very useful and simple change of variables is to replace variables 
by monomials in new variables. {\bf Please note that in what follows, we will 
sometimes use \underline{real} exponents.} 
\begin{dfn}
For any ring $R$, let $R^{m\times n}$ denote the set of $m\times n$ matrices
with entries in $R$. For any $M\!=\![m_{ij}]\!\in\!\R^{n\times n}$ 
and $y\!=\!(y_1,\ldots,y_n)$, we define the formal expression 
$y^M\!:=\!(y^{m_{1,1}}_1\cdots y^{m_{n,1}}_n,\ldots, 
y^{m_{1,n}}_1\cdots y^{m_{n,n}}_n)$. We call the substitution 
$x\!:=\!y^M$ a {\bf monomial change of variables}. Also, for 
any $z\!:=\!(z_1,\ldots,z_n)$, we let $xz\!:=\!(x_1z_1,\ldots,x_nz_n)$. 
Finally, let $\gln(\Z)$ \linebreak \scalebox{.94}[1]{denote the set of all 
matrices in $\Z^{n\times n}$ with determinant $\pm 1$ (the set of {\bf 
unimodular} matrices). \dia } 
\end{dfn} 

\begin{prop} 
\label{prop:monochange} 
(See, e.g., \cite[Prop.\ 2]{tri}.) 
For any $U,V\!\in\!\R^{n\times n}$, we have the formal identity 
$(xy)^{UV}\!=\!(x^U)^V(y^U)^V$. Also, if 
$\det U\!\neq\!0$, then the function $e_U(x)\!:=\!x^U$ is an analytic 
automorphism of $\Rn_+$, and preserves smooth points and singular points of 
zero sets of analytic functions. Finally, 
$U\!\in\!\gln(\Z)$ implies that $e^{-1}_U(\Rn_+)\!=\!\Rn_+$ and that $e_U$ 
maps distinct open orthants of $\Rn$ to distinct open orthants of $\Rn$. \qed 
\end{prop} 

Via a simple application of {\bf Hermite factorization} (see 
Definition \ref{dfn:hermite} and Lemma \ref{lemma:unimod} below), we 
can derive the following corollary which reveals why the dimension 
related hypotheses of Theorem \ref{thm:me} are mild and necessary. 
\begin{cor} 
\label{cor:hyp} 
Given any $n$-variate $m$-nomial $f$ with $d\!=\!\dim \supp(f)\!<\!n$, 
we can find (within $\pp$) a $U\!\in\!\gln(\Z)$ such that 
$g(y)\!:=\!f(y^U)$ is a $d$-variate $m$-nomial 
with $\dim \supp(f)\!=\!d$, and $g$ 
vanishes in $(\Rs)^{d}$ iff $f$ vanishes in $\Rsn$. Moreover, there is an 
absolute constant $c$ such that $\size(g)\!=\!O(\size(f)^c)$. \qed 
\end{cor}
\begin{dfn}
\label{dfn:hermite}
\cite{unimod,storjo}
Given any $M\!\in\!\Z^{m\times n}$, the {\bf Hermite factorization}
of $M$ is an identity of the form $UM\!=\!H$ where $U\!\in\!\glm(\Z)$ and
$H\!=\![h_{ij}]\!\in\!\Z^{n\times n}$ is nonnegative and upper triangular,
with all off-diagonal entries smaller than the positive diagonal entry in
the same column. Finally, the {\bf Smith factorization of $M$} is
an identity
of the form $UMV\!=\!S$ \linebreak 
\scalebox{.98}[1]{where $U\!\in\!\glm(\Z)$, $V\!\in\!\gln(\Z)$, and
$S\!=\![s_{ij}]\!\in\!\Z^{m\times n}$ is diagonal, with 
$s_{i,i}|s_{i+1,i+1}$ for all $i$. \dia}
\end{dfn}
\begin{lemma}
\label{lemma:unimod}
\scalebox{.95}[1]{
\cite{unimod,storjo}
For any $M\!=\![m_{i,j}]\!\in\!\Z^{n\times n}$, the Hermite and Smith
factorizations of $M$}\linebreak exist uniquely,
and can be computed within $O(n^4
\log^3(n\max_{i,j}|m_{i,j}|))$ bit operations.\linebreak 
Furthermore, in the notation of
Definition \ref{dfn:hermite}, the
entries of $U$, $V$, $S$, and $H$ all have bit size
$O(n^3\log^2(2n+\max_{i,j}|m_{i,j}|))$. \qed
\end{lemma}

To prove Theorems \ref{thm:me} and \ref{thm:ckt}, we will first need some 
tricks for efficiently 
deciding when polynomials in $\cA$ have roots in $\Rsn$ or $\Rn_+$.  
\begin{lemma} 
\label{lemma:signs} 
Suppose $f$ is an $n$-variate $(n+1)$-nomial with affinely independent 
support $A\!=\!\{a_1,\ldots,a_{n+1}\}\!\subset\!\Rn$, with 
$a_j\!=\!(a_{1,j},\ldots,a_{n,j})$ for all $j$. Then   
\begin{enumerate} 
\item{$f$ has a root in $\Rn_+$ iff not all the coefficients of 
$f$ have the same sign. }  
\item{If $A\!\subset\!\Zn$ and all the coefficients 
of $f$ have the same sign, then $f$ has a root in 
$\Rsn$ iff there are indices $i,j,j'$ with $a_{i,j}-a_{i,j'}$ odd 
iff there are indices $i,j$ with $a_{i,j}-a_{i,1}$ odd. } 
\end{enumerate} 
\end{lemma} 

\noindent 
While this last lemma is ultimately elementary, we were unable to 
find any similar explicit statement in the literature. So we supply a proof 
in Section \ref{sec:loose}. Another tool we will need is a description 
of certain real zero sets ``at infinity''. 
\begin{dfn}
\label{dfn:face}
Given any compact $S\!\subset\!\Rn$ and $w\!\in\!\Rn$, the {\bf face
of $S$ with inner normal $w$} --- denoted $S^w$ --- is the
set of all $x\!\in\!S$ that minimize the inner product $x\cdot w$.
In particular, a {\bf facet} of $S$ is a face $S^w$ for some
$w$ with $\dim S^w\!=\!\dim S -1$. \dia
\end{dfn}
\begin{dfn}
The {\bf Newton polytope} of $f$, $\newt(f)$, is the convex hull of
the support of $f$. Also, for any $w\!=\!(w_1,\ldots,w_n)\!\in\!\Rn$,
the {\bf initial term function} of
$f(x)\!=\!\sum_{a\in A}c_ax^a$ {\bf with respect to the weight $w$} is
$\init_w(f)\!:=\!\!\!\!\!\!\!\!\sum\limits_{\substack{a\in A\\
w\cdot a \text{ minimized}}} c_ax^a$. \dia
\end{dfn}
\begin{thm}
\label{thm:non}
\scalebox{.95}[1]{\cite[Theorem 3 and Lemmata 14 and 15]{tri}
Let $f$ be any $n$-variate $m$-nomial $f$}\linebreak with
$n$-dimensional Newton polytope and define 
$\cN(f)$ to be the number of connected \linebreak components of $f$.  
Then
\begin{enumerate} 
\item{If $Z_+(\init_w(f))$ is smooth for all $w\!\in\!\Rn\setminus\{\bO\}$, 
then $\cNn(f) \ \leq \ \sum\limits_{\substack{w 
\text{ a unit inner facet } \\
\text{ normal of } \newt(f)}} \!\!\!\!\!\!\!\!\!\!\!\!
\cN(\init_w(f))$.} 
\item{ [$\cN(\init_w(f))\!>\!0$ and 
$Z_+(\init_w(f))$ smooth] for some $w\!\in\!\Rn\setminus\{\bO\} 
\Longrightarrow Z_+(f)$ has a non-compact connected component.} 
\item{ [$f\!\in\!\Z[x_1,\ldots,x_n]$ and $w\!=\!(0,\ldots,0,1)$] 
$\Longrightarrow$ for $\eps\!>\!0$ sufficiently 
small, $Z_+(f)$ has no more than $\cN(\init_w(f)-\eps)$ 
non-compact connected components with limit points in 
$\R^{n-1}_+\times\{0\}$. } 
\item{$\cNn(f)\!\geq\!1 \Longrightarrow \cN(\init_w(f))\!\geq\!1$ 
for some $w\!\in\!\Rn\setminus\{\bO\}$. \qed} 
\end{enumerate} 
\end{thm} 

\noindent
The initial term functions above are also known as {\bf initial term
polynomials} or {\bf face polynomials} when the exponents are integral.

\section{The Proofs of Our Main Results: Theorems 3, 2, and 1}  
\label{sec:proofs}
We will use $[k]$ in place of $\{1,\ldots,k\}$ throughout. Let us 
start by highlighting two 
of the most important theoretical tools behind our proofs: the binomial 
formula for circuit discriminants and Baker's famous result on approximating 
linear forms in logarithms. 

First recall that $A\!\subset\!\Rn$ is a {\bf circuit}
iff $A$ is affinely dependent, but every proper subset of $A$ is affinely
independent.\footnote{ This terminology comes from matroid theory and
has nothing to do with circuits from complexity theory.}
Also, we say that $A$ is a {\bf degenerate circuit} iff
$A$ contains a point $a$ and a proper subset $B$ such that 
$a\!\in\!B$, $A\setminus a$ is affinely independent, and $B$ is a circuit.
For instance, \epsfig{file=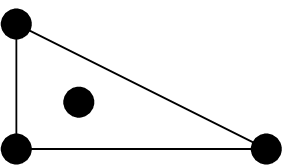,height=.35cm}
and \epsfig{file=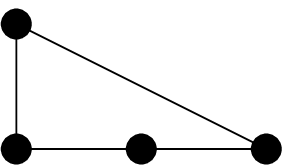,height=.35cm} are
respectively a circuit and a degenerate circuit.

We can now summarize what we 
need to know about $A$-discriminants when $A$ is a circuit:
\begin{lemma}
\label{lemma:ckt} 
Suppose $A\!=\!\{a_1,\ldots,a_{n+2}\}\!\subset\!\Zn$ is a circuit of 
cardinality $n+2$.  Then, defining $\size(A)$ to be the sum of the sizes of 
the coordinates of $A$, and letting $Z^*_\R(f)$ (resp.\ 
$Z^*_\C(f)$) denote the zero set of $f$ in $\Rsn$ (resp.\ 
$\Csn$), we have: 
\begin{enumerate} 
\item{There is a unique (up to sign) vector 
$m\!=\!(m_1,\ldots,m_{n+2})\!\in\!\Z^{n+2}$ such that (a) the coordinates of 
$m$ are all nonzero and have gcd $1$, (b) 
$\sum^{n+2}_{j=1} m_ja_j\!=\!\bO$, \linebreak (c) 
$\sum^{n+2}_{j=1} m_j\!=\!0$, and 
(d) $|m_j|\!:=\!\frac{\vol(\conv(A\setminus\{a_j\}))}{\gcd\{
\vol(\conv(A\setminus\{a_i\}))\; | \; i\in [n+2]\}}$ for all $j$. } 
\item{The vector $m$ from Assertion (1) can be found in time 
polynomial in $\size(A)$. }  
\item{Permuting (if necessary) the $a_j$ so that $m_j\!>\!0$ for all 
$j\!\in\![k]$ and $m_j\!<\!0$ for all $j\!\in\!\{k+1,\ldots,n+2\}$, we 
have\\ 
\mbox{}\hfill$\Delta_A(c_1,\ldots,c_{n+2})\!=\!  
\left(\prod\limits^k_{i=1} m^{m_i}_i\right) 
\left(\prod\limits^{n+2}_{i=k+1} c^{-m_i}_i\right) - 
\left(\prod\limits^{n+2}_{i=k+1} m^{-m_i}_i\right) \left(\prod\limits^k_{i=1} c^{m_i}_i\right)$.\hfill \mbox{} } 
\item{ $\Delta_A(c'_1,\ldots,c'_{n+2})\!=\!0$ for 
some $c'\!=\!(c'_1,\ldots,c'_{n+2})\!\in\!(\Cs)^{n+2} \Longleftrightarrow 
Z^*_\C(f)$ contains a degenerate point $\zeta$ (and 
any such $\zeta$ satisfies $\zeta^{a_i}\!=\!\left.\frac{\partial 
\Delta_A}{\partial c_i}\right|_{c=c'}$ for all $i$). In particular, 
$Z^*_\R(f)$ has at most one degenerate point per open orthant of $\Rsn$. } 
\end{enumerate} 
\end{lemma} 

\noindent 
{\bf Proof:} Assertions (1) and (3) of Lemma \ref{lemma:ckt}
follow immediately from \cite[Prop.\ 1.8, Pg.\ 274]{gkz94}.

Assertion (2) follows
easily from Lemma \ref{lemma:unimod}, upon observing that $m$ is merely
the generator of the integral kernel of a suitable integral matrix. 

\scalebox{.92}[1]{To prove Assertion (4), observe that by Assertion (3) and 
Proposition \ref{prop:monochange}, $\Delta_A(c_1,\ldots,
c_{n+2})\!=\!0$}\linebreak 
for some $c\!:=\!(c_1,\ldots,c_{n+2})\!\in\!(\Cs)^{n+2} \Longrightarrow c$ is 
a smooth point on a hypersurface defined by a real binomial equation.
Via \cite[Thm.\ 1.5, Ch.\ 1, pp.\ 16]{gkz94}, we then obtain that  
$Z^*_\C\left(\sum^{n+2}_{j=1}c_jx^{a_j}\right)$ has a degeneracy 
$\zeta$, and any such $\zeta$ must satisfy the binomial 
system stated above. (In particular, the coordinates are 
{\bf rational} in $c_1,\ldots,c_{n+2}$ if $A$ generates 
$\Zn$ as a lattice, thanks to Lemma \ref{lemma:unimod}.) Moreover, 
by Proposition \ref{prop:monochange} and Lemma \ref{lemma:unimod}, 
we easily derive that a real binomial system can have no more than $1$ solution 
per open orthant of $\Rsn$. 

Conversely, $Z^*_\R\left(\sum^{n+2}_{j=1}c_jx^{a_j}\right)$ has a 
degenerate point $\Longrightarrow \Delta_A(c_1,\ldots,c_{n+2})\!=\!0$, by the 
definition of $\Delta_A$. \qed 
\begin{baker}
\scalebox{.8}[1]{{\bf (Special Case)} [Bak77] 
Suppose $c_0,\ldots,c_N\!\in\!\Z$, $\alpha_1,\ldots,\alpha_N\!\in\!\N$,
$C\!:=\!\max\{4,c_0,\ldots,c_N\}$,}\linebreak and 
$A\!:=\!\log^{N}\max\{4,\alpha_1,
\ldots,\alpha_N\}$. Also let $\Lambda\!:=\!c_0+
c_1\log(\alpha_1)+\cdots+c_N\log(\alpha_N)$. Then
$\Lambda\!\neq\!0\Longrightarrow |\Lambda|\!>\!(CA)^{-(16N)^{200N}
A\log A}$. \qed
\end{baker}

\subsection{The Proof of Theorem \ref{thm:disc}}  
\label{sec:disc} 
Theorem \ref{thm:disc} will follow easily from the 
two algorithms we state immediately below, once we prove their correctness 
and verify their polynomial-time complexity. However, we will first need 
to recall the concept of a gcd-free basis. In essence, a gcd-free 
basis is nearly as powerful as factorization into primes, but is 
{\bf far} easier to compute.  
\begin{dfn}
\label{dfn:gcd}
\cite[Sec.\ 8.4]{bs}
For any subset $\{\alpha_1,\ldots,\alpha_N\}\!\subset\!\N$, a {\bf gcd-free
basis} is a pair of sets $\left(\{\gamma_i\}^\eta_{i=1},
\{e_{ij}\}_{(i,j)\in [N]\times [\eta]}\right)$ such that
(1) $\gcd(\gamma_i,\gamma_j)\!=\!1$ for all $i\!\neq\!j$, and
(2) $\alpha_i\!=\!\prod^{\eta}_{j=1} \gamma^{e_{ij}}$ for all $i$. \dia
\end{dfn}
\begin{thm}
\label{thm:gcd}
\cite[Cor.\ 4.8.2 and Thm.\ 4.8.7 of Sec.\ 4.8]{bs}
Following the notation of \linebreak
Definition \ref{dfn:gcd},
there is a gcd-free basis for $\{\alpha_1,\ldots,\alpha_N\}$, with 
$\eta$, $\size(\gamma_i)$, and $\size(e_{ij})$
each polynomial in $\sum^N_{\ell=1}\size(\alpha_\ell)$, 
for all $i$ and $j$. Moreover, one can always find such a
gcd-free basis using just
$O\!\left(\left(\sum^N_{\ell=1}\size(\alpha_\ell)\right)^2\right)$
bit operations. \qed
\end{thm}
\begin{algbinv}\mbox{}\\
{\bf Input:} Integers $\alpha_1,\beta_1,u_1,v_1,\ldots,\alpha_N,\beta_N,
u_N,v_N$. \\
{\bf Output:} A true declaration as to whether
$\alpha^{u_1}_1\cdots\alpha^{u_N}_N\!=\!\beta^{v_1}_1\cdots\beta^{v_N}_N$. \\
{\bf Description:}\\
\vspace{-.8cm}
\begin{enumerate}
\addtocounter{enumi}{-1}
\item{If $\prod^N_{i=1}(\sign \; \alpha_i)^{u_i \ \mod \ 2}\neq
\prod^N_{i=1}(\sign \; \beta_i)^{v_i \ \mod \ 2}$ then
output {\tt ``They are not equal.''} and stop. Otherwise,
replace the $\alpha_i$ and $\beta_i$ by their absolute values. }
\item{Construct, via Theorem \ref{thm:gcd}, a gcd-free basis
$(\{\gamma_i\}^\eta_{i=1},\{e_{ij}\}_{(i,j)\in [2N]\times [\eta]})$ for
$\alpha_1,\ldots,\alpha_N$,\linebreak
$\beta_1,\ldots,\beta_N$. }
\item{If $\sum^N_{i=1} e_{ij}u_i\!=\!\sum^{2N}_{i=N+1}e_{ij}v_i$ for all
$j\!\in\![\eta]$ then output {\tt ``They are equal.''} and stop.}
\item{Output {\tt ``They are not equal.''} and stop.}
\end{enumerate}
\end{algbinv}
\begin{algbins}\mbox{}\\
{\bf Input:} Positive integers $\alpha_1,\beta_1,u_1,v_1,\ldots,\alpha_N,
\beta_N,u_N,v_N$.\\
{\bf Output:} The sign of $\alpha^{u_1}_1\cdots\alpha^{u_N}_N -
\beta^{v_1}_1\cdots\beta^{v_N}_N$. \\
{\bf Description:}\\
\vspace{-.8cm}
\begin{enumerate}
\item{Let
$C\!:=\!\log^{2N}\max\{4,\alpha_1,\beta_1,\ldots,\alpha_N,\beta_N\}$,
$M\!:=\!\max\{4,u_1,v_1,\ldots,u_N,v_N\}$, and \linebreak
$E\!:=\!\frac{1}{3}(CM)^{-(32N)^{400N}C\log C}$. }
\item{For all $i\!\in\![N]$, let $A_i$ (resp.\ $B_i$) be an approximation
of $\log\alpha_i$ (resp.\ $\log \beta_i$) within $\frac{E}{2NM}$
(using, say, Arithmetic-Geometric Mean Iteration \cite{dan}).\footnote{
Other approximation techniques can be used as well: It is sufficient to
use any algorithm that can find the $b$ leading bits of $\log N$ within
a number of bit operations polynomial in $b+\log N$. } }
\item{Output the sign of
$\left(\sum^N_{i=1}u_iA_i\right)-\left(\sum^N_{i=1}v_iB_i\right)$ and stop.}
\end{enumerate}
\end{algbins}

\noindent
Note that while we can certainly compute $A\!:=\!\alpha^{u_1}_1\cdots
\alpha^{u_N}_N$ using a number of {\bf arithmetic} operations polynomial
in $s\!:=\!\size(u_1)+\cdots+\size(u_N)$, the {\bf bit} size of $A$
is already exponential in $s$; hence the need for
our last two algorithms.

\begin{lemma}
\label{lemma:bin} 
Algorithms {\tt BinomialVanish} and {\tt BinomialSign} are both correct. 
Moreover, Algorithm {\tt BinomialVanish} runs in time polynomial 
in\\ 
\mbox{}\hfill 
$I\!:=\!\sum^N_{i=1}(\size(\alpha_i)+\size(\beta_i)+\size(u_i)+\size(v_i))$ 
\hfill \mbox{}\\
and, if $N$ is fixed, Algorithm {\tt BinomialSign} runs in time polynomial 
in $I$ as well. 
\end{lemma}

\noindent 
{\bf Proof of Lemma \ref{lemma:bin}:}  
That Algorithm {\tt BinomialVanish} is correct and runs in
time polynomial in $I$ follows directly from Theorem \ref{thm:gcd}. 

That Algorithm {\tt BinomialSign} is correct and runs in time 
polynomial in $I$ (for {\bf fixed} $N$) follows easily 
from Baker's Theorem: First, taking logarithms, observe 
that the sign of $\alpha^{u_1}_1\cdots\alpha^{u_N}_N-\beta^{v_1}_1\cdots 
\beta^{v_N}_N$ is the same as the sign of 
$S\!:=\!\left(\sum^N_{i=1}u_i\log\alpha_i\right)- 
\left(\sum^N_{i=1}v_i\log\beta_i\right)$. 
Clearly then, $\left|S-\left[\left(\sum^N_{i=1}u_iA_i\right)-
\left(\sum^N_{i=1}v_iB_i\right)\right]\right|\!<\!E$, so 
Baker's Theorem tells us that Step (3) of Algorithm {\tt BinomialSign} 
indeed computes the sign of $S$.  So we have correctness. 

To see that Algorithm {\tt BinomialSign} runs in time polynomial in $I$ 
for fixed $N$, first note that the $A_i$ and $B_i$ each require 
$O(\log(M)+\log(N)-\log E)
\!=\!O\left((32N)^{400N}C\log (C)\log(CM)\right)$ bits of accuracy.\footnote{  
Note also that the true number of bits of accuracy we would use in 
practice is $2^\mu$, where $\mu$ is the smallest integer with, 
say, $2^\mu\!>\!(32N)^{400N}C\log(C)\log(CM)$. This is because while 
it may be non-trivial to compute $E$ exactly, finding $2^\mu$ is easy via the 
old trick of recursive squaring.}  
So, via our chosen method for approximating logarithms \cite{dan}, 
we see\linebreak \scalebox{.85}[1]{that the complexity of our algorithm is 
polynomial in 
$O\left(N(32N)^{400N}C\log (C)\log(CM)+\sum^N_{i=1}\log(\alpha_i\beta_i)
\right)$}\linebreak
$=\!2^{O(N\log N)}\left(\sum^N_{i=1}\size(\alpha_i)+\size(\beta_i)
\right)^{O(N)} \sum^N_{i=1}(\size(u_i)+\size(v_i))$, and we are done. \qed 

\medskip
\noindent
{\bf The Proof of Theorem \ref{thm:disc}:} First note 
that if our input $A$ is not a circuit,\footnote{One can in fact check in 
polynomial-time whether $A$ is a circuit. See, e.g., Proposition \ref{prop:pp} 
of Section \ref{sub:proof1}.} then $A$ is a degenerate circuit, and 
$\Delta_A(f)$ is then identically $1$. This is because $Z^*_\R(f)$ is smooth 
when the support of $f$ is a degenerate circuit 
(see, e.g., the proof of Case 2 of Theorem \ref{thm:ckt} in Section 
\ref{sub:proof2}). 
So we can assume that $A\!=\!\{a_1,\ldots,a_{n+2}\}$ is a circuit, $f$ 
has support $A$, and that $c_j$ is the coefficient of $x^{a_j}$ in $f$ for 
all $j$. 

Thanks to Lemma \ref{lemma:ckt}, Assertion (1) (resp.\ 
Assertion (2)) follows straightforwardly from the complexity bound for  
Algorithm {\tt BinomialVanish} (resp.\ {\tt BinomialSign}) 
we just proved in Lemma \ref{lemma:bin}. 
In particular, the latter lemma tells us that the complexity of 
$\disc$, for an input $(A,c_1,\ldots,c_{n+2})$, is  
polynomial in $\sum^{n+2}_{i=1}\log(c_im_i)$ (following the notation of 
Lemma \ref{lemma:ckt}); and the same is true for 
$\sdisc$ provided $n$ is fixed. The classical Hadamard matrix 
inequality \cite{mignotte} tells us that  
$\size(m_i)\!=\!O(n\log(n\max_{j,k}\{a_{jk}\}))$, so 
the complexity of $\disc$ is indeed polynomial 
in $\size(f)$; and the same holds for $\sdisc$ when $n$ is fixed. \qed 

\subsection{Deforming to Polyhedra: The Proof of Theorem \ref{thm:ckt}} 
\label{sub:proof2} 
Let $A\!=\!\{a_1,\ldots,a_m\}$ be the support of $f$ and 
and write $f(x)\!=\!\sum^m_{j=1}c_jx^{a_j}$. Since $Z_+(f)$ is 
unaffected if $f$ is replaced by a monomial multiple of $f$, we can clearly 
assume without loss of generality that $a_1\!=\!\bO$. 

If $n\!=\!1$ then $f$ must be a univariate trinomial and 
Theorem \ref{thm:ckt} follows immediately from Descartes' Rule. So 
let us assume henceforth that $n\!\geq\!2$.

\medskip 
\noindent
\scalebox{.93}[1]{\fbox{{\bf Case 1: $A$ affinely independent:}}   
Letting $A'$ denote the matrix whose columns are $a_2,\ldots,a_m$,}\linebreak 
observe that $\barf(y)\!:=\!f(y^{A'^{-1}})\!=\!c_1+c_2y_1+\cdots+c_my_m$,
with $m\!\leq\!n$. Moreover, thanks to Proposition \ref{prop:monochange},
$Z_+(\barf)$ (the intersection of a hyperplane with the 
positive orthant of $\Rn$) and $Z_+(f)$ are diffeomorphic, so we are done. In 
particular, we see that $Z_+(f)$ is either empty or a connected, open,  
$C^\infty$, real $(n-1)$-manifold. \qed
\begin{rems} 
\scalebox{.93}[1]{Note that we can in fact allow arbitrary {\bf real} exponents 
for $f$ in our proof above. \dia} 
\end{rems} 

\noindent
\scalebox{.91}[1]{\fbox{{\bf Case 2: $A$ is a degenerate circuit:}}   
Suppose, without loss of generality, that $B\!=\!\{\bO,a_2,\ldots,
a_\ell\}$}\linebreak (with $\ell\!<\!m\!\leq\!n+2$) is a circuit, 
and $a_\ell$ in the relative interior of $\conv B$ if $B$ intersects 
the relative interior of $\conv B$. 

Let $A''$ be the $n\times (m-2)$ 
matrix whose columns are $a_{m-1},\ldots,a_2$ and (via Lemma \ref{lemma:unimod})
define $U$ to be any $n\times n$ unimodular matrix $U$ such that 
$UA''$ is lower triangular, with nonnegative diagonal.  
Defining $\bar{f}(y)\!:=\!f(y^U)$, 
it is then easily checked that $\bar{f}(y)$ is of the form $c_1+c_2y^{a''_2}
+\cdots+ c_my^{a''_m}$ with $a''_j\!\in\!\{0\}^{n+1-\max\{j,m-2\}}\times
\Z^{\max\{j,m-2\}-1}$ for all $j\!\in\!\{2,\ldots,m-1\}$, and 
$a''_m\!\in\!\{0\}^{n+2-m}\times (\Z\setminus\!\{0\})\times \Z^{m-3}$. 
Moreover, Proposition \ref{prop:monochange} tells us that
$Z_+(f)$ and $Z_+(\bar{f})\times \R^{n-m}_+$ (when   
$Z_+(\barf)$ is considered as a subset of $\R^m_+$) are diffeomorphic.

Defining $(a''_{1,j},\ldots,a''_{n,j})\!:=a''_j$ for all $j$ and
$\bar{f}_1(y)\!=\!\frac{-(c_1+c_2y^{a''_2}+\cdots+c_{m-1}y^{a''_{m-1}})}
{y^{a''_{m,1}}_1 \cdots y^{a''_{m,m-1}}_{m-1}}$, observe then\linebreak 
\scalebox{.97}[1]{that $\bar{f}_1\!\in\!\Z[y_1,\ldots,y_{m-1}]$ and, for all 
$y\!=\!(y_1,\ldots,y_m)\!\in\!\R^m_+$, we have 
$\bar{f}(y)\!=\!0 \Longleftrightarrow
\bar{f}^{1/a''_{m,m}}_1\!=\!y_m$.}\linebreak So
we see that $Z_+(\bar{f})$ is exactly the positive part of the graph
of the $C^\infty$ function $\bar{f}^{1/a_{m,n}}_1 : \R^{m-1}_+\longrightarrow
\R$. Clearly then, $Z_+(\barf)$ consists of a union of
connected, open, $C^\infty$, real $(m-1)$-manifolds, 
and $Z_+(f)\!\approx\!Z_+(\barf)\times\R^{n-m}_+$ is thus smooth and has 
{\bf no} compact connected components. So let us now bound 
$\cN_{\mathrm{non}}(\barf)$.  

If $\barf_1$ is always positive on $\R^{m-1}_+$ 
then $Z_+(\barf)$ clearly consists of a single non-compact 
connected component. So we can henceforth assume that $Z_+(\barf)$ is 
non-empty, which in turn implies that {\bf every} connected component $C$ of 
$Z_+(\barf)$ has a limit point in $\R^{m-1}_+\times\{0\}$. Assertion (3) of 
Theorem \ref{thm:non} then tells us that $\cN_{\mathrm{non}}(\barf)\!\leq\! 
\cN(\barf_1-\eps)$ for some $\eps\!>\!0$. So, by induction (reducing to 
lower-dimensional instances of Cases 2 or 3), we are done. \qed 

At this point, we must recall a result of Viro on the classification 
of certain real algebraic hypersurfaces. In what follows, we liberally 
paraphrase from \cite[Thm.\ 5.6]{gkz94}. 
\begin{dfn} 
Given any finite point set $A\!\subset\!\Rn$, let us call any function 
$\omega : A \longrightarrow \R$ a {\bf lifting}, 
denote by $\pi : \R^{n+1}\longrightarrow \Rn$ the natural 
projection which forgets the last coordinate, and let 
$\hat{A}\!:=\!\{(a,\omega(a))\; | \; a\!\in\!A\}$. We then say that the 
polyhedral subdivision $\Sigma_\omega$ of $A$ defined by 
$\{\pi(Q) \; | \; Q \text{ a 
lower$^{10}$ face of } \conv \hat{A} \text{ of dimension } 
\dim A\}$ is {\bf induced by the lifting $\omega$}, and we call 
$\Sigma_\omega$ a {\bf triangulation induced by a lifting} iff every cell of 
$\Sigma_\omega$ is a simplex.\footnotetext{A {\bf lower} face is simply a 
face which has an inner normal with positive last coordinate.} Finally, 
given any $f(x)\!=\!\sum_{a\in A}c_ax^a\!\in\!\Z[x_1,\ldots,x_n]$, we 
define $f_{\omega,\eps}(x)\!:=\!\sum_{a\in A}c_a\eps^{\omega(a)}x^a$ to be 
the {\bf toric perturbation of $f$ (corresponding to the lifting $\omega$)}. 
\dia 
\end{dfn} 
\begin{dfn} 
Following the notation above, suppose $\dim A\!=\!n$ and $A$ 
is equipped with a triangulation $\Sigma$ induced by a lifting {\bf and} a 
function $s : A \longrightarrow \{\pm\}$ which we will call a 
{\bf distribution of signs for $A$}. We then define a locally 
piece-wise linear manifold  --- the {\bf Viro diagram} $\cV(\Sigma,s)$ --- in 
the following local manner: For any $n$-cell $C\!\in\!\Sigma$, 
let $L_C$ be the convex hull of the set of midpoints of edges of 
$C$ with vertices of opposite sign, and then define 
$\cV(\Sigma,s)\!:=\!\bigcup\limits_{C \text{ an } n\text{-cell}} 
L_C$. \dia 
\end{dfn} 
\begin{ex} 
The following figure illustrates $6$ circuits of cardinality $4$, 
each equipped with a triangulation induced by a lifting, and a distribution 
of signs. The corresponding (possibly empty) Viro diagrams are drawn 
in the lightest color visible (yellow on the color version of this 
paper). \dia\\ 
\mbox{}\hfill\epsfig{file=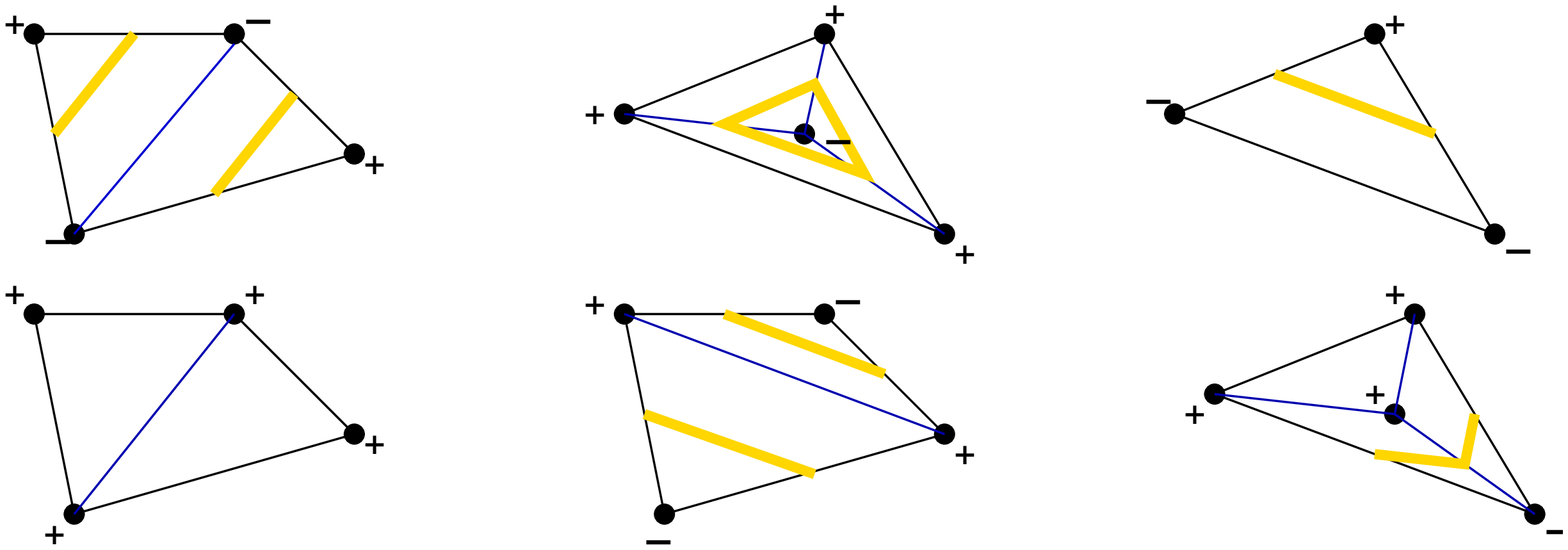,height=2.3in}\hfill 
\end{ex} 
\begin{thm} 
\label{thm:viro} 
Suppose $f(x)\!=\!\sum_{a\in A}c_ax^a\!\in\!\Z[x_1,\ldots,x_n]$ with 
$\supp(f)\!=\!A$ and $\dim A\!=\!n$, $\omega$ is any lifting of $A$, and 
define $s_f(a)\!=\!\sign(c_a)$ for all $a\!\in\!A$. Then for any sufficiently 
small $\eps\!>\!0$, $Z_+(f_{\omega,\eps})$ is isotopic to 
$\cV(\Sigma_\omega,s_f)\setminus\partial \conv A$. In particular, 
$\cV(\Sigma_\omega,s_f)$ is a disjoint finite union of piece-wise linear 
manifolds, each possibly having a non-empty boundary. \qed 
\end{thm} 
\begin{lemma} 
\label{lemma:low} 
Suppose $A$ is a circuit, $\Sigma$ is a triangulation of $A$, $n\!=\!\dim A$, 
and $s$ is any distribution of signs on $A$. Call any point $a$ in the 
relative interior of $A$ with $s(a)$ opposite $s(a')$ for all $a'\!\in\!A
\setminus\{a\}$ a {\bf caged alternation of $(A,s)$}. Then\\ 
1. $Z_+(f)$ smooth $\Longrightarrow Z_+(f)$ is isotopic to 
$\cV(\Sigma,s_f)\setminus\conv A$ for some $\Sigma$. \\
2. $\cV(\Sigma,s)$ has no boundary iff $\cV(\Sigma,s)$ is the boundary of 
an $n$-simplex iff [$A$ has a caged\linebreak 
\mbox{}\hspace{.6cm}alternation $a$ and $\Sigma$ is the 
triangulation obtained by the lifting that sends 
$a\mapsto 0$ and\linebreak
\mbox{}\hspace{.6cm}$a'\mapsto 1$ for all $a'\in A\setminus\{a\}$]. 
\end{lemma} 

\noindent 
{\bf Proof of Lemma \ref{lemma:low}:} By Lemma \ref{lemma:ckt} and 
Lemma \ref{lemma:unimod}, it easily follows that $A$ 
has at most $2$ discriminant chambers in $\R^{n+2}_+$, and each such 
chamber contains a unique toric perturbation. Since the topology of $Z_+(f)$ 
is constant on any discriminant chamber containing $f$ (e.g., \cite[Ch.\ 11, 
Sec.\ 5A, Prop.\ 5.2, pg.\ 382]{gkz94}), we obtain Assertion (1). 
 
Now note that by definition, any triangulation of a circuit $A$ obtained by 
lifting has at most $n+1$ top-dimensional cells (since it is the projected 
lower hull of a $(n+1)$-simplex). So then, if $\cV(\Sigma,s)$ has no 
boundary, its convex hull 
must clearly have dimension $n$, in which case we see that $\cV(\Sigma,s)$ 
is a union of at least $n+1$ simplices of dimension $n-1$. So 
$\cV(\Sigma,s)$ is the union of exactly $n+1$ simplices of dimension 
$n-1$, equal to the boundary of its convex hull, i.e., 
$\cV(\Sigma,s)$ is the boundary of a $n$-simplex. This proves 
the first rightward implication of Assertion (2), and the converse is obvious. 

Now if $\cV(\Sigma,s)$ is the boundary of a $n$-simplex, then 
$A$ clearly intersects its relative interior, which in turn 
implies that $\Sigma$ must be the specified triangulation. 
(Indeed, it is a standard fact that any circuit has exactly 
$2$ triangulations. So the only other possible $\Sigma$ 
has exactly $1$ $n$-cell and could not possibly give the 
$\cV(\Sigma,s)$ we desire.) Furthermore, $A$ must then 
clearly contain a caged alternation, for otherwise $\cV(\Sigma,s)$ 
would no longer be the boundary of an $n$-simplex. This proves 
the second rightward implication of Assertion (2), and the converse is 
obvious. \qed 

\noindent 
We are now ready to return to our proof of Theorem \ref{thm:ckt} and 
finish the remaining special case. 

\medskip 
\noindent
\fbox{{\bf Case 3: $A$ is a circuit:}} First note that by 
Corollary \ref{cor:hyp}, we can apply a monomial change of 
variables and assume $m\!=\!n+2$ without loss of generality. 
(Moreover, if $m\!<\!n+2$ initially, then every connected 
component of $Z_+(f)$ must be non-compact, thanks to Proposition 
\ref{prop:monochange}.) 

Let us first bound $\cNn(f)$: First note that 
since $A$ is a circuit, every facet of $\conv(A)$ has affinely independent 
vertex set. So by Case 1, the facet functions of $f$ each automatically
have smooth zero sets which are either empty or consist of a
single non-compact connected component. The Upper
Bound Theorem of polyhedral combinatorics \cite{edelsbrunner} then implies
that $\conv(A)$ has no more facets than a moment $n$-polytope with $n+2$
vertices. The latter polytope has exactly $4$, $6$, or $9$ facets, 
according as $n$ is $2$, $3$, or $4$ \cite{edelsbrunner}. 
So Assertion (1) of Theorem \ref{thm:non} then directly implies our stated 
bound on $\cNn(f)$, for $n\!\in\!\{3,4\}$. That 
$\cNn(f)\!\leq\!2$ for $n\!=\!2$ follows from earlier work of 
Daniel Perrucci using a different argument, involving a 
detailed analysis of the central special case $1+x+y+Ax^ay^b$ 
\cite[Thm.\ 4, Assertion (4)]{perrucci}.  

Note in particular that when $n\!=\!1$, $Z_+(f)$ has no non-compact components 
unless $f$ is identically zero. Note also that $f(x,y)\!:=\!(x-1)(y-1)\!=\!
xy-x-y+1$ has support a circuit and exactly two connected components 
(each non-compact) for its zero set in $\R^2_+$. So our stated bounds 
for $\cNn(f)$ are indeed tight for $n\!\in\!\{1,2\}$. 

Assume now that $n\!\geq\!5$. 
If $Z_+(f)$ is smooth, then Assertion (1) of Lemma \ref{lemma:low} tells us 
that $Z_+(f)$ is isotopic to some Viro diagram. As observed within the 
proof of Lemma \ref{lemma:low}, a triangulation $\Sigma$ of $A$ contains 
no more than $d+1$ $d$-simplices, and thus a Viro diagram 
of the form $\cV(\Sigma,s)$ conists of no more than $d+1$ $(d-1)$-simplices. 
Since any such $(d-1)$-simplex can belong to at most one connected 
component of $\cV(\Sigma,s)$, and since $\cV(\Sigma,s)$ is a union 
of piece-wise linear manifolds with boundary, we see that $\cV(\Sigma,s)
\setminus\conv A$ has at most $d+1$ non-compact connected components. 
So by Theorem \ref{thm:viro}, $\cNn(f)\!\leq\!n+1$ in the smooth case. 

Now recall the following two standard inequalities:\\
($N_c$) \hfill $\cNc(f)\!\leq\!\cNc(f-\eps)+\cNc(f+\eps)$
\hfill\mbox{}\\
($N_n$) \hfill $\cNn(f)\!\leq\!\cNn(f-\eps)+\cNn(f+\eps)$
\hfill\mbox{}\\
for $\eps\!>\!0$ sufficiently small. Moreover,
$Z_+(f\pm\eps)$ is smooth for all $\eps\!>\!0$ sufficiently small.
(See, e.g., \cite[Lemma 2]{basu}.)

Having proved $\cNn(f)\!\leq\!n+1$ in the smooth case, Inequality ($N_n$) 
then immediately implies that $\cNn(f)\!\leq\!2n+2$, even in the presence 
of singularities for $Z_+(f)$. 

To bound $\cNc(f)$, let $n$ be arbitrary once again. Observe then that Lemma 
\ref{lemma:low} implies that $\cNc(f)\!\leq\!1$, with equality 
only if $Z_+(f)$ is isotopic to an $(n-1)$-sphere, as long as  
$Z_+(f)$ is smooth. So our bound for $\cNc(f)$ holds if 
$Z_+(f)$ is smooth. 

On the other hand, if $Z_+(f)$ has a singularity, 
consider first the special case where $(A,s_f)$ does {\bf not} have a 
caged alternation. Lemma \ref{lemma:low} and ($N_c$) then imply that 
$\cNc(f)\!=\!0$. So we have in fact proved a 
strengthening of the final assertion of Theorem \ref{thm:ckt}. 

As for the case where $(A,s_f)$ has a caged alternation,  
assume without loss of generality that it is $a_{n+2}$. 
Lemma \ref{lemma:low} then tells us that $\cNc(f)\!\leq\!1$ 
(with $Z_+(f)$ isotopic to an $(n-1)$-sphere if $\cNc(f)\!=\!1$), 
provided $Z_+(f)$ is smooth. So we have our assertion for compact components 
of $Z_+(f)$ in the smooth case. 
 
To conclude, assume that $Z_+(f)$ has a singularity. 
Lemma \ref{lemma:ckt} then tells us that this singularity is 
unique and that $\left(\prod^{n+1}_{i=1}
m^{m_i}_i\right)c^{-m_{n+2}}_{n+2}\!=\!m^{-m_{n+2}}\prod^{n+1}_{i=1}
c^{m_i}_i$. So $f-\eps$ and $f+\eps$ thus lie in {\bf opposite} 
discriminant chambers for any $\eps\!>\!0$. (Note that we are still 
assuming that $(A,s_f)$ has a caged alternation.) So one of $Z_+(f-\eps)$ 
or $Z_+(f+\eps)$ is empty. By Lemma \ref{lemma:low},  
and inequalities ($N_c$) and ($N_n$), we thus obtain that 
$\cNc(f)\!\leq\!1$ and $\cNn(f)\!=\!0$. 
In particular, $Z_+(f)$ must have consisted of a single point, for 
otherwise, both $Z_+(f-\eps)$ and $Z_(f+\eps)$ would have been 
non-empty (by the Implicit Function Theorem). 

Having proved our upper bound for $\cNc(f)$, we need only 
exhibit examples proving tightness for each $n\!\geq\!1$. For 
$n\!=\!1$ there is the obvious example of $(x-1)^2\!=\!x^2-2x+1$. 
As for $n\!\geq\!2$, Theorem \ref{thm:viro} tells us that 
it suffices to use $\eps(1+x^{2n}_1+\cdots+x^{2n}_n)-x_1\cdots x_n$, 
for any $\eps\!>\!0$ sufficiently small. So we are done. \qed

\subsection{Phase Transitions: The Proof of Theorem \ref{thm:me}}
\label{sub:proof1}
The complexity lower bounds of Theorem \ref{thm:me}  --- Assertions 
(3) and (4) --- are the easiest to prove, so we start there: 
                
\noindent 
\fbox{{\bf The Proof of Assertion (3):}} 
Recall that $\sat$ is the problem of deciding whether a 
Boolean formula of the form $B(X)=C_1(X)\wedge \cdots \wedge C_k(X)$ has a 
satisfying assignment, where $C_i$ is of one of the following forms:\\
\mbox{}\hfill $X_i\vee X_j \vee X_k$, \   
$\neg X_i\vee X_j \vee X_k$, \  
$\neg X_i\vee \neg X_j \vee X_k$,  \ 
$\neg X_i\vee \neg X_j \vee \neg X_k$, \hfill \mbox{}\\ 
$i,j,k\!\in\![3N]$, and a satisfying assigment consists of 
an assigment of values from $\{0,1\}$ to the variables 
$X_1,\ldots,X_{3N}$ which makes the equality $B(X)\!=\!1$ true. 
$\sat$ is one of the most basic $\np$-complete problems 
\cite{gj}. In particular, for our purposes, let us 
measure the size of a $\sat$ instance such as the one above as $N$.

Let us now observe that $\feas_\R(\cM)$ is $\np$-hard, where 
$\cM$ is the family of all polynomial systems of the form 
$(f_1,\ldots,f_k)$ where, for all $i$, 
$f_i\!\in\!\Z[x_1,\ldots,x_{3N}]$ involves no more 
than $3$ variables and has degree $\leq\!1$ with respect to each of 
them.\footnote{While 
this construction is well-known in some circles, we have included 
a detailed explanation for the convenience of the reader.} 
To see why, first observe that any Boolean formula $B(X)$ with 
{\bf no} occurence of ``$\wedge$'' can be 
converted into a polynomial $f_B(x)\!\in\!\Z[x_1,\ldots,x_{3N}]$ via the 
following table of substitutions:

\noindent 
\mbox{}\hfill\begin{tabular}{ccc} 
$X_i$ & $\mapsto$ &  $x_i$ \\
$\neg X_i$ & $\mapsto$ & $1-x_i$ \\
$C_i(X)\vee C_j(X)$ & $\mapsto$ & 
$f_{C_i}(x)+f_{C_j}(x)-f_{C_i}(x)f_{C_j}(x)$ 
\end{tabular} \hfill\mbox{}\\

\noindent 
For instance, the formula $X_1\wedge \neg X_3\wedge X_8$ 
becomes\\
\mbox{}\hfill $(x_1+(1-x_3)-x_1(1-x_3))+x_8-(x_1+(1-x_3)-x_1(1-x_3))x_8\!=\!
1-x_3+x_1x_3+x_3x_8-x_1x_3x_8$.\hfill\mbox{}\\ 
To any $\sat$ instance as above,  
we can then associate the polynomial system\\
\mbox{}\hfill $(f_{C_1}-1,\ldots,f_{C_k}-1,
x_1(1-x_1),\ldots,x_{3N}(1-x_{3N}))$,\hfill\mbox{}\\ 
and it is easily checked that $B$ has a satisfying assignment iff  
$F_B$ has a real root. (Moreover, any root of $F_B$ clearly 
lies in $\{0,1\}^{3N}$.) Note also that the size of $F_B$ 
is clearly $O(N)$, and that $F_B$ has at least as many equations as 
variables. Clearly then, $\feas_\R(\cM)\!\in\!\pp 
\Longrightarrow \sat\!\in\!\pp$, and thus $\feas_\R(\cM)$ is 
$\np$-hard.  

\scalebox{.94}[1]{To conclude, we need only recall that any system of 
polynomials $F$ chosen from $\Z[x_1,\ldots,x_n]$}\linebreak 
can always be converted to its 
{\bf Shor Normal Form} $S_F$ \cite{shor}. In particular, 
$F$ has a real root iff $S_F$ has a real root, the size of 
$S_F$ is linear in the size of $F$, $S_F$ consists of 
linear trinomials and/or binomials of degree $\leq\!2$, and 
the number of new equations introduced is the same as the number of 
new variables introduced. (More concretely, 
substitutions like $x_2-x^7_1  \mapsto (y_1-x^2_1,y_2-y^2_1,y_3-y_2y_1,
x_2-y_3x_1)$ can be used to reduce all powers to $2$ or less, and 
substitutions like $y_k\!=\!c_ix^{a_i}+c_jx^{a_j}$ can be used to reduce 
any polynomial to a collection of polynomials, each with $3$ or fewer 
monomial terms.)  
So the number of variables of $S_F$ is bounded above by 
the number of equations of $S_F$, and we are done. \qed 

\medskip 
\noindent 
\fbox{{\bf The Proof of Assertion (4):}} By our proof of Assertion (3), and 
replacing any system $(f_1,\ldots,f_k)\!\in\!\cS$ with the polynomial 
$f^2_1+\cdots+f^2_k$, it is easy to see that $\feas_\R(\cE)$ is $\np$-hard, 
where $\cE$ is the family of $n$-variate $11n$-nomials that are sums of 
squares. Indeed, by our earlier reduction to $\sat$, we can additionally 
assume that any polynomial in $\cE\cap\Z[x_1,\ldots,x_n]$ contains  
$x^2_1(1-x_1)^2,\ldots,x^2_n(1-x_n)^2$ as summands. The latter 
assumption easily implies that $\dim \supp(f)\!=\!n$ for any 
$n$-variate polynomial $f\!\in\!\cE$, for then $\newt(f)$ 
must contain a line segment parallel to each and every 
coordinate axis of $\Rn$. 

To reduce to the family $\cC_\eps$, simply note that from any 
$f\!\in\!\cE$ and $\eps\!>\!0$, we can form the new polynomial 
$g_\eps(x,y)\!:=\!f(x)+y^2_1+\cdots+y^2_N$, where $N\!:=\!\lceil 10n^{1/\eps}
\rceil$.  
Observe then that $g$ involves exactly $n+N$ variables, 
$g$ has exactly $11n+N$ monomial terms, and $\dim \supp(g)\!=\!n+N$. 
In particular,\\ 
\mbox{}\hfill $11n+N - (n+N) \!=\!10n\!\leq\!\lceil 10n^{1/\eps}\rceil^\eps
\!\leq\!(n+N)^\eps$,\hfill\mbox{}\\ 
for all $n\!\geq\!1$. Moreover, for any fixed $\eps\!>\!0$, the size of 
$g$ is clearly polynomial in the size of $f$. So we then clearly obtain that 
$\feas_\R(\cC_\eps)\!\in\!\pp \Longrightarrow \feas_\R(\cE)\!\in\!\pp$, and 
we are done. \qed 

\medskip 
The complexity upper bounds of Theorem \ref{thm:me} --- Assertions (1) and 
(2) --- then follow easily from the two final algorithms we state below. 
In what follows, we let $Z_\R(f)$ denote the real zero set of $f$ and 
$A\!:=\!\supp(f)$. 
\begin{algsimp}\mbox{}\\
{\bf Input:} A polynomial $f\!\in\!\Z[x_1,\ldots,x_n]$ with 
affinely independent support $A\!=\!\{a_1,\ldots,a_m\}$ of 
cardinality $m$, and $c_j$ the coefficient of $x^{a_j}$ in 
$f$ for all $j$. \\
{\bf Output:} True declarations as to whether $Z_+(f)$, 
$Z^*_\R(f)$, and $Z_\R(f)$ are, respectively, empty or not. \\
{\bf Description:}\\
\vspace{-.8cm} 
\begin{enumerate} 
\item{If the $c_i$ are {\bf not} all of the same sign, 
then output\\ 
\mbox{}\hfill
{\tt ``$Z_+(f)$, $Z^*_\R(f)$, and $Z_\R(f)$ are all 
non-empty.''}\hfill\mbox{}\\ and stop. } 
\item{If any one of $a_2-a_1,\ldots,a_m-a_1$ has an 
odd coordinate, then output\\ 
{\tt ``$Z^*_\R(f)$ and $Z_\R(f)$ 
are non-empty, but $Z_+(f)$ is empty.''}\hfill\mbox{}\\ 
and stop. } 
\item{Output {\tt ``$Z_+(f)$ and $Z^*_\R(f)$ are empty.''}. 
\item{If $\bO\!\not\in\!A$ then output {\tt ``$Z_\R(f)$ is non-empty.''} 
and stop.} 
\item{Output {\tt ``$Z_\R(f)$ is empty.''} and stop. }  
 \qed}  
\end{enumerate} 
\end{algsimp}
\begin{rems}
\label{rems:nc}
Note that $\dim \supp(f)$ can be computed in $\nc_1$ as follows:
compute the rank of the matrix whose columns are $a_2-a_1,
\ldots,a_m-a_1$ via the parallel algorithm of Csanky \cite{csanky}. So checking whether a given $f$ is a valid input to Algorithm
{\tt DimPlus1Case} can be done within $\nc_1$. \dia
\end{rems}

\noindent
\fbox{{\bf The Proof of Assertion (1):}} Here we simply apply Algorithm 
{\tt DimPlus1Case}. Assuming Algorithm {\tt DimPlus1Case} is correct, the
complexity upper bound then follows trivially, since one
can check all the necessary signs and/or parities in
logarithmic time, using a number of processors linear in $\size(f)$. 
(Note also that
we can check within $\nc_1$ whether an input $f$ has $\supp(f)$ affinely
independent, thanks to Remark \ref{rems:nc}.) So we need only show that
Algorithm {\tt DimPlus1Case} is correct.

To prove the latter, note that Lemma \ref{lemma:signs} implies that
the existence of a pair of coefficients of $f$
with opposite sign is the same as $Z_+(f)$ being non-empty.
So Step (1) is correct, and we can assume henceforth that all
the coefficients of $f$ have the same sign.
By Lemma \ref{lemma:signs} again, the
existence of indices $i,j$ with $a_{i,j}-a_{i,1}$ odd implies
that\linebreak $Z^*_\R(f)\setminus Z_+(f)$ is non-empty.
So Step (2) is correct, assuming its hypothesis is true.

So let us now assume that $a_{i,j}$ has the same parity as
$a_{i,1}$ for all $i,j$. Note that if $f$ has a root in $\Rn$ then this root
must lie in some coordinate subspace $L$ of
minimal positive dimension. So, by our initial hypotheses, on $L$, the
polynomial $f$ will restrict to an $n'$-variate $m'$-nomial with
$m'\!\leq\!n'+1$ and affinely independent support a subset of $A$. In
particular, since we now have that all the coefficients of $f$ have the same
sign, and since $a_{i,j}$ has the same parity as $a_{i,1}$ for all $i,j$,
one final application of Lemma \ref{lemma:signs} implies that $Z_\R(f)$ is
empty. So Step (2) is also correct when its hypothesis is false.

Steps (3)--(4) are clearly correct since $f$ vanishes at $\bO$ iff 
$f$ does {\bf not} have a nonzero constant term. \qed
\begin{rems}
Let $\feas^+_\R$ and $\feas^*_\R$ denote the obvious analogues of
$\feas_\R$ where we respectively restrict to roots in $\Rn_+$ or
$\Rsn$. Also, paralleling our earlier notation, let $\feas^+_\R(\cF)$
and $\feas^*_\R(\cF)$ be the corresponding natural restrictions of 
$\feas^+_\R$ and $\feas^*_\R$ to inputs in some family $\cF$.
Our preceding proof then clearly implies that $\feas^+_\R(\cA)\!\in\!\nc_1$
(even if real exponents are allowed) and $\feas^*_\R(\cA)\!\in\!\nc_1$. \dia
\end{rems}

We can now describe our final algorithm:  
\begin{algcirc}\mbox{}\\ 
{\bf Input:} A polynomial $f\!\in\!\Z[x_1,\ldots,x_n]$ with
support $A\!=\!\{a_1,\ldots,a_{n+2}\}$ of cardinality $n+2$, 
$\dim A\!=\!n$, $a_1\!=\!\bO$, and an $\ell\!\in\!\{3,\ldots,n+2\}$ 
such that $B\!:=\!\{a_1,\ldots,a_\ell\}$ is a circuit and $a_\ell$ 
lies in the relative interior of $\conv B$ 
if $B$ intersects the relative interior of $\conv B$. \\
{\bf Output:} True declarations as to whether $Z_+(f)$ 
has compact and/or non-compact connected components, 
along with notification as to whether $Z_+(f)$ is a point. 
\\
{\bf Description:}\\
\vspace{-.6cm} 
\begin{enumerate}
\addtocounter{enumi}{-1}
\item{Let $W$ denote the set of all inner facet normals, with 
integer coordinates having no common factor, of $A$. Also let $c_j$ be the 
coefficient of $x^{a_j}$
in $f$ for all $j$, and let $A''$ be the $n\times n$ matrix whose columns 
are $a_{n+1},\ldots,a_2$.}  
\item{If $n\!\geq\!2$, decide (using Algorithm {\tt DimPlus1Case}) whether 
there is a $w\!\in\!W$ with $A^w$ affinely independent and $Z_+(\init_w(f))$ 
non-empty. If so, then output\\ \scalebox{.98}[1]{{\tt ``$Z_+(f)$ is non-empty 
and all its connected components are non-compact.''}}\linebreak and stop. }
\item{If $A$ is a degenerate circuit then do the following:  
\begin{enumerate}  
\item{If $B$ intersects the relative interior of $\conv B$ then do the 
following:
\begin{enumerate} 
\item{ Find (via Lemma \ref{lemma:unimod} of Section \ref{sec:back}) a 
unimodular matrix $U$ such that $UA''$ is lower triangular and has a 
nonnegative diagonal with exactly one zero entry. Then, replacing $f$ by $-f$ 
if necessary, assume that all the $c_j$ (except possibly $c_\ell$) are 
positive. Finally, define $h(x)\!:=\!f(x)-c_{\ell+1}x^{a_{\ell+1}}-\cdots 
-c_{n+2} x^{a_{n+2}}$ and $g(y)\!:=\!h(y^U)$. } 
\item{Decide, via a lower-dimensional instance of 
Algorithm {\tt DimPlus2Case}, whether $Z_+(g)$ contains at least 
$2$ points. If so, then then output\\ 
\scalebox{.8}[1]{{\tt ``$Z_+(f)$ is non-empty, smooth, and all its 
connected components are non-compact.''}}\linebreak and stop. }   
\end{enumerate}} 
\end{enumerate}} 
\item{If $A$ intersects the interior of $\conv A$ then let $s_f(a)\!:=\! 
\sign(c_a)$ for all $a\!\in\!A$, order the sequence of signs 
$s_f(a)$ in increasing order of $a$ (if $n\!=\!1$), and 
do the following:
\begin{enumerate} 
\item{If $n\!=\!1$ and $s_f$ has exactly one 
sign alternation then output\\ 
{\tt ``$Z_+(f)$ is a point of multiplicity $1$.''} and stop.} 
\item{ If $(A,s_f)$ has a caged alternation then do the following: 
\begin{enumerate} 
\item{
Decide, via Algorithm {\tt 
BinomialVanish}, whether $\Delta_A(c_1,\ldots,c_{n+2})\!=\!0$. 
If so, then output\\ 
\scalebox{.85}[1]{{\tt ``$Z_+(f)$ has exactly one connected 
component, and it is a singular point.''}}\\ and stop. }
\item{Decide, via Algorithm {\tt BinomialSign}, whether 
$\sign(\Delta(g))\!=\!(-1)^{m_{n+2}}$. If so, then 
output\\ {\tt ``$Z_+(f)$ has exactly one connected component, 
and it is smooth}\linebreak {\tt and isotopic to an $(n-1)$-sphere.''} 
and stop. } 
\end{enumerate} }
\end{enumerate} } 
\item{Output {\tt ``$Z_+(f)$ is empty.''} and stop. \qed } 
\end{enumerate}
\end{algcirc}

It will be useful to observe that the input hypotheses to 
the preceding algorithm can be checked within $\pp$: 
\begin{prop}
\label{prop:pp}
Given any finite set $A\!\subset\!\Zn$ with $\#A\!\leq\!\dim A+2$, we can
decide if $A$ contains a circuit $B$, 
and find the unique such $B$ should it exist, in 
time polynomial in $\size(A)$. 
Moreover, if $B$ contains a circuit and $B$ intersects the 
relative interior of $\conv B$, then we can find the 
unique point of this intersection also in time polynomial in $\size(A)$. 
\end{prop} 

\noindent 
{\bf Proof of Proposition \ref{prop:pp}:} 
First, one simply checks via the method of Remark
\ref{rems:nc} if $\# A\!=\!\dim A+2$ (for if not,
$A$ can not contain a circuit). Then, via Cramer's Rule and
the Newton identities \cite[Ch.\ 15, Pgs.\ 292--296]{bcss},
one simply checks which subsets of $A$ of cardinality $\dim A$
form facets of $A$. (Overlaps can be distinguished via a computation 
of the underlying inner facet normals.) If all the facets of $A$
have cardinality $\dim A$, then $A$ itself is a circuit
and we are done. Otherwise, some facet $S$ of $A$
has cardinality $>\!\dim A$, and we use the same
method recursively to find the unique circuit of $S$.
Via \cite[Prop.\ 21, Ch.\ 15, pp.\ 295]{bcss}, it is then easily checked 
that the bit complexity of this method is no worse than $O(\size(A)^8)$, 
where $\size(A)$ is the sum of the sizes of the coordinates of the points of 
$A$.  

To efficiently check whether $B$ intersects the relative interior of
$\conv B$, we can employ any linear
programming algorithm with polynomial-time bit complexity\footnote{
See, e.g., \cite[Ch.\ 15]{bcss} for a nice description of a barrier
method employing Newton's method.} as follows:
express --- whenever possible --- each
$b_i\!\in\!B$ as a convex linear combination of points in $B\setminus\{b_i\}$.
If some point $b_i$ can be expressed in this way, then this
point is unique, and we can permute the entries of $A$ so that 
$a_1\!=\!\bO$, 
$a_\ell\!=\!b_i$, and $a_j\!=\!b_j$ for all $j\!\in\![\ell]\setminus\{i\}$.
Otherwise, $B$ does not intersect its relative interior. \qed 

\medskip 
\noindent
\fbox{{\bf The Proof of Assertion (2):}} 
First note that $f$ has a nonzero constant term iff $f$ does {\bf not} have
$\bO$ as a root, and this can be checked with just $1$ bit operation.
{\bf So we can assume henceforth $f$ has a nonzero constant term.}
Also note that the polynomial obtained from $f$ by setting any subset of its
variables to $0$ lies in $\cB_{n'}\cup \cA$ for some $n'\!<\!n$.
Moreover, since we can apply changes of variables like
$x_i\mapsto -x_i$ in $\pp$, and since there are exactly $3^n$
sequences of the form $(\eps_1,\ldots,\eps_n)$ with
$\eps_i\!\in\!\{0,\pm 1\}$ for all $i$, it suffices at this point to show 
that $\feas^+_\R(\cB_n)\!\in\!\pp$ for {\bf fixed} $n$.

This will be accomplished by Algorithm {\tt DimPlus2Case}, and 
via Proposition \ref{prop:pp} and a monomial change of variables 
(employing Lemma \ref{lemma:unimod}), we can indeed prepare $f$ to be a 
suitable input to this algorithm.\footnote{We should note that 
if $\dim\supp(f)\!<\!n$ initially then Corollary \ref{cor:hyp} 
implies that every connected component of $Z_+(f)$ will be non-compact.} 
So we can now assume that $f$ satisfies 
these input hypotheses. Assertion (2) then follows immediately --- assuming 
that Algorithm {\tt DimPlus2Case} is correct and runs in $\pp$ 
for fixed $n$.  So let us now prove correctness and analyze the 
complexity along the way. 

\noindent
{\bf \underline{Correctness and Complexity of Steps (0)--(1):}} 
The correctness of Steps (0)--(1) follows immediately from Theorem 
\ref{thm:non}. Note also that $W$ consists of no more than 
$\begin{pmatrix} n+2 \\ n\end{pmatrix}\!=\!(n+2)(n+1)/2$ normals, 
and each such normal can be constructed (employing Cramer's Rule and 
the Newton identities \cite[Ch.\ 15, Pgs.\ 292--296]{bcss}) 
via $n-1$ determinants of $(n-1)\times (n-1)$ matrices, followed by  
a gcd computation (see, e.g., \cite{bs}[Ch.\ 3] for a detailed exposition 
on near optimal gcd algorithms). Since 
we've already proved that Algorithm {\tt DimPlus1Case} runs in $\pp$, 
Steps (0)--(1) clearly run in $\pp$, and no further complexity analysis is 
needed if the hypothesis of Step (1) is satisfied.  

If the hypothesis of Step (1) is false, then (applying Lemma \ref{lemma:signs} 
again) this means that 
all the $c_j$ with $a_j$ a vertex of $\conv A$ have the same sign. 
So we can assume without loss of generality that the $c_j$ are all positive  
and continue to Step (2). \qed 

\noindent
{\bf \underline{Correctness and Complexity of Step (2):}} 
\scalebox{.9}[1]{First let us define 
$h_j(x)\!:=\!f(x)-\sum^{n+2}_{i=n+3-j}c_{i}x^{a_i}$}\linebreak   
and $g_j(y)\!:=\!h_j(y^U)$ for all $j\!\in\![n+2-\ell]$. Note that   
$g_j\!\in\!\R[x_1,\ldots,x_{n-j}]$ for all $j\!\in\![n+2-\ell]$. 

By the Case 2 portion of the Proof of Theorem \ref{thm:ckt}, it immediately 
follows that $Z_+(f)$ is smooth and diffeomorphic to the positive part of the 
graph of $-g_1$ (as a function on $\R^{n-1}_+$). Clearly then, 
every connected component of $Z_+(f)$ is non-compact. Furthermore, 
since $c_1,\ldots,c_{\ell-1}\!>\!0$ (and $c_\ell\!>\!0$ as well, if 
$a_\ell$ is a vertex of $\conv A$), Lemma \ref{lemma:signs} and 
Assertion (2) of Theorem \ref{thm:non} imply that every connected component 
must have a limit point on $\R^{n-1}\times\{0\}$. 
Moreover, it follows easily from the Implicit Function Theorem and Assertion 
(4) of Lemma \ref{lemma:ckt} that $Z_+(g_1)$ has at least $2$ points 
$\Longrightarrow Z_+(f)$ is non-empty. So by induction, 
$Z_+(f)$ is non-empty iff $Z_+(g_{n+1-\ell})$ is non-empty (provided 
$\ell\!<\!n+1$). 

Note also that by construction, $\supp(g_j)$ is a degenerate 
circuit for all $j\!\in\![n+1-\ell]$ (provided $\ell\!<\!n+1$). So, to 
simplify notation, we can clearly assume without loss of generality that 
$\ell\!=\!n+1$ and that $\supp(g)$ is a circuit. Just as 
in the last paragraph, we can still assert the two implications 
(a) $Z_+(g)$ has at least 
$2$ points $\Longrightarrow Z_+(f)$ is non-empty, and (b) 
$Z_+(g)$ empty $\Longrightarrow Z_+(f)$ is empty. So we are left  
with the special case where $Z_+(g)$ is a point. Moreover, 
unless $(B,s_{g})$ has a caged alternation, the coefficients of 
$g$ will all be positive, thus making $Z_+(f)$ empty (and this 
will be correctly declared later in Step (4)). 
So we can assume that $Z_+(g)$ is a point {\bf and} $(B,s_g)$ 
has a caged alternation. 

Define $g_\eps(y)\!:=\!\eps\left(\prod^{\ell-1}_{i=1}c_iy^{a''_i}\right)
+c_{\ell}y^{a''_{\ell}}$. Assertion (3) of Lemma \ref{lemma:ckt} then 
immediately implies that for $\delta,\eps\!>\!0$ sufficiently small, 
$\sign(\Delta_B(g-\delta))\!=\!\sign(\Delta_B(g_\eps))$.  
In particular, $g-\delta$ must then lie in the same $B$-discriminant 
chamber as $g_\eps$, and Theorem \ref{thm:viro} then implies that 
the graph of $-g$ attains a maximum value of $0$ within $\R^{n-1}_+$.  
In other words, $Z_+(f)$ is empty, and this would be correctly 
declared later in Step (4). So Step (2) 
is correct, and its complexity is dominated by an instance 
in $\cB_{n'}$ of Algorithm {\tt DimPlus2Case}, for some $n'\!<\!n$. 

We can thus assume now that $A$ is a circuit and continue to 
Step (3). \qed  

\noindent
{\bf \underline{Steps (3) and (4):}} 
The correctness of Step (3-a) follows immediately from Descartes' 
Rule of Signs, combined with the observation that the derivative 
of $f(x)$ (or $x^{a_3}f(1/x)$) is nonzero on $\R_+$. 

If the hypothesis of Step (3-b) is violated, then all the coefficients 
of $f$ must have the same sign and $Z_+(f)$ must then be empty 
(and this will be correctly declared later in Step (4)). 
So we may assume that $(A,s_f)$ has a caged alternation. Moreover, 
we may also assume without 
loss of generality that $c_1,\ldots,c_{n+1}\!>\!0$ and $c_{n+2}\!<\!0$. 

If $Z_+(f)$ is smooth then Assertion (1) of Lemma \ref{lemma:low} implies 
that  $Z_+(f)$ is isotopic to one of two possible Viro diagrams, easily 
seen to be either empty or the boundary of an $n$-simplex. 
Moreover, by Theorem \ref{thm:viro}, $Z_+(f)$ is isotopic to the 
latter diagram iff $\sign(\Delta_A(f))\!=\!\sign(\Delta_A(f_\eps))$ 
for all $\eps\!>\!0$ sufficiently small, 
where $f_\eps(x)\!:=\!\eps\left(\prod^{n+1}_{i=1}c_ix^{a_i}\right)
+c_{n+2}x^{a_{n+2}}$. A simple calculation from Assertion (3) of Lemma 
\ref{lemma:ckt} then tells us that $\sign(\Delta_A(f_\eps))\!=\!(-1)^{m_{n+2}}$ 
for all $\eps\!>\!0$. In other words, Step (3-b-ii) is correct. 

So now assume $Z_+(f)$ has a singularity $\zeta$. Then, 
Assertion (4) of Lemma \ref{lemma:ckt} implies that $\zeta$ is the {\bf only} 
singularity of $Z_+(f)$. Moreover, since $(A,s_f)$ has a 
caged alternation, $Z_+(f)$ must then be exactly $\{\zeta\}$, as 
already proved toward the end of the Case 3 portion of the proof of 
Theorem \ref{thm:ckt}. So Step (3-b-i) is correct. 

To conclude, observe that any case not satisfying the hypotheses of 
any of our steps results in $Z_+(f)$ being empty, and this is 
correctly declared by Step (4). Furthermore, we see that 
the complexity of Steps (3)--(4) is dominated by a single instance 
of Algorithm {\tt BinomialVanish} and a single instance of 
{\tt BinomialSign}, for input $f$. \qed 

\section{The Proof of Lemma \ref{lemma:signs}} 
\label{sec:loose} 
As before, we can assume without loss of 
generality that $a_1\!=\!\bO$ and $c_1\!=\!1$ by dividing by a suitable 
monomial. Furthermore, we can clearly permute the $a_i$ so that 
$c_i\!>\!0$ iff $i\!\leq\!k$, for some $k\!\leq\!n+1$. 
Let $A'$ be the matrix whose columns are $a_2,\ldots,a_{n+1}$. 
Then, via the change of variables 
$x\!\!=\!z/(|c_2|,\ldots,|c_{n+1}|)^{A^{-1}}$ (which clearly preserves 
the existence of roots of $f$ in any open orthant of $\Rsn$), 
we can then clearly assume that $f(x)\!=\!1+x^{a_2}+\cdots+x^{a_k}-
x^{a_{k+1}}- \cdots -x^{a_{n+1}}$. 

Our criteria for checking the existence of roots of $f$ in $\Rn_+$ 
or $\Rsn$ then clearly reduce to checking whether $k\!<\!n$ or 
whether $A'$ has an odd entry. So let us prove that the latter 
conditions correctly characterize the existence of roots of $f$ in 
$\Rn_+$ and $\Rsn$. 

First note that $f(x)\!=\!0$ for some $x\!\in\!\Rsn$ iff\\ 
\mbox{}\hfill($\star$) \ \ \ \ \ \ \ $x^{A'}\!=\!\alpha \text{ and } 
1+\alpha_1+\cdots+\alpha_k-\alpha_{k+1}-\cdots-\alpha_n\!=\!0$,\hfill
\mbox{} \\
for some $\alpha\!=\!(\alpha_1,\ldots,\alpha_n)\!\in\!\Rsn$.  
Assertion (1) then follows almost trivially: Assuming $x\!\in\!\Rn_+$, 
the equality $k\!=\!n$ and Proposition \ref{prop:monochange} imply that 
$f(x)\!=\!1+\alpha_1+\cdots +\alpha_n\!>\!0$, 
so there can be no roots for $f$ in $\Rn_+$. Taking the inverse implication, 
suppose $k\!<\!n$. Then we can set $\alpha\!:=\!
\left(\underset{k}{\underbrace{1,\ldots,1}},\underset{n-k}
{\underbrace{\frac{k+1}{n-k},\ldots,\frac{k+1}{n-k}}}\right)$ to obtain 
$1+\alpha_1+\cdots+\alpha_k-\alpha_{k+1}-\cdots-\alpha_n\!=\!0$. 
So if we can solve $x^{A'}\!=\!\alpha$ over $\Rn_+$, we will have found a root 
in $\Rn_+$ for $f$. Proposition \ref{prop:monochange} tells us that we can 
indeed (since $\det A'\!\neq\!0$), so we are done.  

We now focus on Assertion (2). Letting $y\!:=\!x^U$, note that 
\[ x^{A'}\!=\!\alpha \Longleftrightarrow 
y^S\!=\!y^{UA'V}\!=\!(x^{A'})^V\!=\!\alpha^V, \] 
thanks to Proposition \ref{prop:monochange}, where $S\!=\![s_{ij}]$ is an 
$n\times n$ diagonal matrix with $s_{1,1}|s_{2,2}|\cdots|s_{n,n}$. So we'll 
be able to find a root in $\Rsn$ for $f$ iff\\ 
\mbox{}\hfill($\heartsuit$) \ \ \ \ $\text{There are } 
\alpha,y\!\in\!\Rsn \text{ with } y^S\!=\!\alpha^V \text{ and } 
1+\alpha_1+\cdots+\alpha_k-\alpha_{k+1}- \cdots -\alpha_n\!=\!0. 
$\hfill\mbox{}\\ 
Let us now separately prove the two directions of the equivalence in  
Assertion (2): 

\noindent 
{\bf ($\Longleftarrow$):} First note that the hypothesis 
is invariant under a common translation of $a_1,\ldots,a_{n+1}$. 
So we can assume $a_1\!=\!\bO$ (and $k\!=\!n$ as given), and our 
hypothesis then translates into $A'\!\subset\!\Zn$ having at least one odd 
entry.  $A'$ having at least one odd entry then implies that 
the mod $2$ reduction of $A'$ has positive $(\Z/2\Z)$-rank, and this 
in turn implies that $s_{1,1}$ is odd. (Since left and right multiplication 
by matrices in $\gln(\Z)$ 
preserves $(\Z/2\Z)$-rank, and $s_{1,1}|s_{2,2}|\cdots |s_{n,n}$.) 
Since the map $e_V(x)\!:=\!x^V$ is clearly an automorphism of the 
open orthants of $\Rsn$ (provided $V\!\in\!\gln(\Z)$), there must then 
clearly be some open orthant (having exactly $j$ positive coordinates) 
which is mapped bijectively onto $\R_-\times\R^{n-1}_+$ 
under $e_V$.  So then define $\alpha$ to be any permutation of the vector 
$\left(\underset{j}{\underbrace{1,\ldots,1}},\underset{n-j}
{\underbrace{-\frac{j+1}{n-j},\ldots,-\frac{j+1}{n-j}}}\right)$ 
such that $\sign(\alpha^V)\!=\!(-1,1,\ldots,1)$. 
Clearly then, $y^S\!=\!\alpha^V$ has a solution in $\Rsn$ and thus, by 
($\heartsuit$) and our choice of $\alpha$, $f$ indeed has a root in $\Rsn$.  

\noindent 
{\bf ($\Longrightarrow$):} We will prove the contrapositive.  
Via translation invariance again, in the notation above, we see that our 
hypothesis is equivalent to all the entries of $A'$ being even, and thus 
all the $s_{i,i}$ must be even. We then obtain, via Proposition 
\ref{prop:monochange}, 
that $y^S\!=\!\alpha^V$ has {\bf no} roots in $\Rsn$ unless 
$\alpha\!\in\!\Rn_+$. But then $\alpha\!\in\!\Rn_+$ implies that 
$1+\alpha_1+\cdots+\alpha_n\!>\!0$ (since $k\!=\!n$ by assumption), so 
there can be {\bf no} roots for $f$ in $\Rsn$. \qed

\section*{Acknowledgements}
The authors thank Francisco Santos for discussions on counting regular
triangulations, and Frank Sottile for pointing out reference \cite{bbs}.  
Thanks also to Dima Pasechnik for discussions, and Sue Geller and Bruce 
Reznick for detailed commentary, on earlier versions of this work. 
The authors also thank the anonymous referees for their suggestions, 
especially the considerable improvement of Assertion (4) of Theorem 
\ref{thm:me}.  

\bibliographystyle{acm}

\end{document}